%% file: main.tex
\documentclass[manuscript,screen]{acmart}

\usepackage[english]{babel}
\usepackage{amsmath}
\usepackage{graphicx}
\usepackage[utf8]{inputenc}
\usepackage{amsmath, mathtools}
\usepackage{algpseudocode}

\usepackage{amsfonts}
\usepackage{amstext}
\usepackage{amsthm}
\usepackage{algorithm}
\usepackage{listings}

\usepackage{subcaption}
\usepackage{graphicx}

\usepackage{pgfplots}
\pgfplotsset{compat=1.18}

\usepackage[utf8]{inputenc}

\usepackage{tikz}

\usepackage{changepage}

\usepackage{comment}

\usepgfplotslibrary{external}
\AtBeginDocument{%
  }

\setcopyright{acmlicensed}
\copyrightyear{2025}
\acmYear{2025}
\acmDOI{XXXXXXX.XXXXXXX}

\begin{document}

\title{Mixed-Precision in High-Order Methods: the Impact of Floating-Point Precision on the ADER-DG Algorithm}

\hyphenation{
Exa-HyPE
}

\author{Marc Marot-Lassauzaie}
\affiliation{%
  \institution{Technical University of Munich}
  \city{Munich}
  \country{Germany}}
\email{marot.marc@tum.de}

\author{Michael Bader}
\affiliation{%
  \institution{Technical University of Munich}
  \city{Munich}
  \country{Germany}}
\email{bader@cit.tum.de}

\renewcommand{\shortauthors}{Marot-Lassauzaie et al.}

\begin{abstract}
We present a mixed-precision implementation of the high-order discontinuous Galerkin method with ADER time stepping (ADER-DG) for solving hyperbolic systems of partial differential equations (PDEs) in the hyperbolic PDE engine ExaHyPE. 
The implementation provides a simple API extension for specifying the numerical precision for individual kernels, and thus allows for testing the effect of low and mixed precision on the accuracy of the solution. 
To showcase this, we study the impact of precision on the overall convergence order and actual accuracy of the method as achieved for four common hyperbolic PDE systems and five relevant scenarios that feature an analytic solution. 
For all scenarios, we also assess how sensitive each kernel of the ADER-DG algorithm is to using double, single or even half precision. 
This addresses the question where thoughtful adoption of mixed precision can mitigate hurtful effects of low precision on the overall simulation.
\end{abstract}


\received{9 April 2025} 
\received[revised]{\dots}  
\received[accepted]{\dots}   

\maketitle

\input{chapters/Introduction}

\input{chapters/SOTA}

\input{chapters/ADERV2}

\input{chapters/Implementation}

\input{chapters/Scenarios}

\input{chapters/Results}

\input{chapters/Conclusion}

\section*{Acknowledgements}\label{chapter:ack}
The presented work was funded by the Deutsche Forschungsgemeinschaft (DFG, German Research Foundation) as part of the project no.\ 462423388, \textit{ExaHyPE MVP -- mixed and variable precision for an Exascale Hyperbolic PDE Engine}.
We thank Tobias Weinzierl and his group, at Durham University, for collaboration, hints and support on the ExaHyPE engine and the Peano framework. 
We also thank Michael Dumbser (University of Trento) for comments and discussions regarding the ADER-DG algorithm. 

\bibliographystyle{ACM-Reference-Format}
\bibliography{sample}

\end{document}

%% file: chapters/Introduction.tex
\section{Introduction}\label{chapter:introduction}

The latest C++23 standard has, for the first time, officially defined so-called \emph{fixed width floating-point types}, providing two ``half precision'' types (stored with only 16 bits) in addition to the existing 32-, 64- and 128-bit formats (see the overview in table \ref{table:precision_bits}). This new introduction illuminates the substantial importance that the deliberate and appropriate choice of floating-point precision has gained for any performance-oriented software.  
In high performance computing (HPC), exploitation of lower precision in research software has been driven by the widening adoption and availability of GPU computing and has been further pushed by the rise in popularity of machine learning applications, where extremely low precision suffices for training and inference \cite{Bader:SIAMnews,Matsuoka:2023}. 
Google for instance has developed Tensor Processing Units able to compute 16- and 8-bit operations specifically to accelerate the inference phase of neural networks \cite{JouppiTPU:2017}. NVIDIA released its own version of Tensor Cores, which perform fused-multiply-add operations on half precision input. 
These specialised low-precision units not only substantially accelerate deep learning applications, but also present an avenue to accelerate other HPC applications \cite{Markidis_tensor_cores}. 

Even on ``standard'' HPC hardware, exploiting low precision promises performance improvement due to lower costs for data transfer (communication and memory), higher cache locality and higher intra-instruction level parallelism (vectorisation, e.g.). To mitigate the obvious downside of low-precision computing, that is the corresponding loss of accuracy,  precision needs to be chosen adaptively in algorithms and implementations.
Such \emph{mixed precision algorithms} are now widely adopted in dense and sparse linear algebra software \cite[e.g.]{Abelfattah_mpc_survey}. A wide range of respective algorithms \cite[e.g.]{higham_mary_2022}, from BLAS routines, LU factorizations, GMRES, etc., up to related methods in machine learning, and multiple libraries facilitate their usage in applications. 
In other fields of scientific computing, however, particularly for the numerical solution of problems posed via (systems of) partial differential equations (PDEs), the adoption of low and mixed precision is still far from being mainstream.

This paper presents a mixed-precision implementation of the ADER-DG algorithm, which combines high-order discontinuous Galerkin (DG) discretisation~\cite{dumbser_ader_2006,dumbser_2008} with \textbf{A}rbitrary high-order \textbf{DER}ivatives (ADER) 
time stepping (introduced by Toro and Titarev~\cite{toro_2002}). 
In particular, we use the ADER-DG algorithm by Dumbser et al.~\cite{dumbser_ader_2014} (with a-posteriori subcell limiting) as implemented in the hyperbolic PDE engine ExaHyPE~\cite{Reinarz_ExaHyPE_2020}. 
ADER-DG has been applied to a wide variety of different problems, including earthquake simulation~\cite{dumbser_ader_2006,wolf_aderdg_seismic}, astrophysics~\cite{Dumbser_CCZ4curl-cleaning2020} or fluid dynamics~\cite{Guerrero:2021,Zanotti_pos_lim}. 
The ExaHyPE engine provides a generic ADER-DG implementation that users can tailor to their particular problem by specifying PDE properties and problem specifics through a simple Python API, facilitating the development of new models. 
Respective applications include the simulation of seismic waves~\cite{Duru:2022}, problems in astrophysics \cite{Zhang_ExaGRyPE2025}, atmospheric flows~\cite{Krenz_PPAM2019} or tsunami propagation~\cite{Leo_thesis,Seelinger:SC21}. 

The presented mixed-precision implementation provides a simple means to configure the precision separately for each major kernel of the ADER-DG algorithm. 
This allows for adapting the applied precision to specific PDEs and scenarios, but also enables users to easily evaluate the impact of lower or higher precision on accuracy in (parts of) the computation, before specifying the precision for production simulations.
We exploit this ease of testing to present a comprehensive study on how low- and mixed-precision computation affects the ADER-DG algorithm, when solving a suite of five scenarios covering four different PDE systems (elastic and acoustic wave equations, Euler equations, shallow water equations). 
We find that numerical precision has a strong impact on $p$-convergence: to achieve high-order convergence on fine meshes, 64-bit precision is a must. 
32-bit precision can suffice for moderately high order and resolution. 
Both 16-bit precision formats were found to be insufficient for many (if not most) scenarios, causing failure of the algorithm in multiple cases. 
Depending on the PDE system and scenario, both lack of mantissa bits (leading to rounding errors) and lack of exponent bits (limiting the range of values) can cause such failure.  
Mixed precision can often ``restore'' such failure cases, even though the low target precision then still dominates the resulting error level.

While performance improvement is a key motivation for using low or mixed precision, we entirely focus on accuracy in this work, and leave performance impact (which can be expected to depend strongly on actual hardware) for future studies. 
Still, we review related work on accuracy and performance impact in chapter~\ref{chapter:state_of_the_art}. 
In chapter~\ref{chapter:ader}, we summarize the key steps of the ADER-DG algorithm. 
The implementation of this algorithm in multiple and mixed precisions is then discussed in chapter~\ref{chapter:implementation}. In chapter \ref{chapter:verification}, we present our suite of evaluation scenarios, stemming from four different hyperbolic PDE systems. 
Chapter~\ref{chapter:results} investigates for these scenarios the impact of computing the entire ADER-DG algorithm in a lower precision. Finally, in chapter \ref{chapter:mixed_precision} we study the impact of using mixed precision, applying lower or higher precision in specific kernels, to identify which ADER-DG kernels are most suited for low precision, or which demand high precision. 

\begin{table}
\centering
\begin{tabular}{|| c | c | c | c c c ||}
 \hline
 Standard & Short & C++ & Mantissa bits & Exponent bits & Max.\ exponent \\ [0.5ex] 
 \hline\hline
 bfloat 16 & bf16 & \texttt{bfloat16} & 7  & 8 & 127 \\
 IEEE binary 16 & fp16 & \texttt{fp16} & 10  & 5 & 15 \\
 IEEE binary 32 & fp32 & \texttt{float} & 23  & 8 & 127 \\
 IEEE binary 64 & fp64 & \texttt{double} & 52  & 11 & 1023 \\
 IEEE binary 128 & fp128 & \texttt{long double} & 112  & 15 & 16383 \\
 \hline
\end{tabular}

\caption{Exponent and mantissa bits in different signed floating-point formats.
In the text, we will also use the colloquial terms ``half'' (for bf16 and fp16), ``single'' (fp32), ``double'' (fp64) and ``quadruple'' (fp128) precision.
}
\label{table:precision_bits}
\end{table}

%% file: chapters/SOTA.tex
\section{Related Work}\label{chapter:state_of_the_art}




Benefits of lowering the precision of floating point values stem from improvements to the hardware's memory usage, decreased latencies in communication, or from more effective usage of the floating point units.
Early general purpose computing on GPUs had to primarily rely on single precision. 
In that context, Göddeke et al.\ performed Finite Element simulations on CPUs and GPUs \cite{GPUFEM}, using both mixed precision (via an iterative refinement approach) and emulation of double precision. They reported speedups of 1.7 for using mixed precision, compared to the double-precision of their code. 
For the first GPU generation to offer half precision, Ho et al.~\cite{Ho_half_precision} developed a code conversion framework to deal with the tedious programming for 16-bit precision. Evaluating a portfolio of mini-apps, they achieved speedups of up to 3 over the single precision version. 
%
%
Recent GPUs offer special instructions, such as for small matrix multiplications, which routinely accelerate computation beyond the 2x speedup normally expected over single-precision. 
In the Turing generation of Nvidia GPUs, tensor cores can offer 4x speedup over the FP16 units and results of higher accuracy \cite{Yan_tensor_cores}. 
Similar benefits can be seen in Yang et al.'s work on exploiting TensorFlow and Cloud TPUs for statistical ferromagnetism simulations using the two-dimensional Ising model \cite{YangMonteCarloIsing}. These simulations were memory bound in all cases, but using half precision both accelerated the code and allowed for larger models, benefiting from reduced memory usage. 
Nowadays, machine learning dominates the hardware development. 
There, the benefits of lower precision are so significant that even 8-bit precision is used (see Wang et al.~\cite{8bitNeural}, e.g.). 

In the context of higher-order Finite Element solvers for elliptic PDEs, Schneck et al.~\cite{schneck_2021} studied the impact of reduced precision storage of preconditioners on the overall error of additive Schwarz smoothers. 
For a test setup solving the Lamé-Navier equation of linear solid mechanics, they found that for a scheme computed with 64 bit compute-precision, 16 bit storage-precision was oftentimes sufficient to achieve convergence (with benefits largest for higher orders). 
Another finding is that for lower precision, fixed-point representations lead to more accurate results than floating point storage. 
Langou et al.~\cite{Langou2006ExploitingTP} studied iterative refinement for linear systems, using a single-precision solver as a kind of preconditioner for a refinement iteration, computed in double precision. 
This method only converges if the condition of the problem does not exceed the reciprocal of the single precision accuracy. 
When tested on a benchmark of 147 matrices 96 converged, 15 did not converge, and 36 failed in one of the calculation steps. 
Carson and Higham extend iterative refinement to a three-precision version for the GMRES algorithm \cite{3PrecGMRES}. 
They show that this method can achieve errors at the working precision, and suggest that on architectures where half precision is implemented efficiently this should perform twice as fast as  single-precision LU factorization. 
This was pushed further by Amestoy et al.~\cite{5PrecGMRES}, who present a five-precision iterative refinement for GMRES. 

With respect to implementation of mixed-precision solvers, Grützmacher et al.\ present a memory-accessor to accelerate memory-bound BLAS operations \cite{AnztGingkoMemoryAccessor}. 
They decouple the precision of arithmetic and memory operations by converting the data to a lower precision format before storing it in memory. 
Data is thus transferred to and from memory faster, yielding higher performance for memory-bound applications while maintaining the high-precision arithmetic. 
%
Demidov et al.~\cite{Demidov_2021} use C++ metaprogramming to adopt blocked and mixed-precision data structures for sparse matrices. For large-scale linear solvers for a Stokes problem, they observed 30\% speedup due to mixed precision (in addition to gains from blocked layouts).  
Current work by Radtke et al.~\cite{PawelMemoryAccessor} suggests memory-centric extensions to C++, which add support for arbitrary floating-point precision for MPI communication. 
Using these extensions yields a speedup of almost two in bandwidth-bound cases, due to halving the memory footprint of the transferred data . 

For hyperbolic PDE systems, mixed-precision approaches are comparably rare, as they typically do not employ solver steps as natural candidates for mixed precision. 
Field et al.~\cite{FieldGPUWENOBlackHole} present a GPU-accelerated mixed-precision WENO method to simulate extremal black holes and gravitational waves. Compared to their base case, which requires quadruple precision, using double precision for all expensive but numerically less critical parts offers a 3.3x speedup without meaningfully impacting the results.
For more general PDE solvers, Croci and de Souza have developed mixed-precision Runge-Kutta-Chebyshev methods \cite{CrociMPCRungeKutta}. These methods employ most of the Runge-Kutta stages for stability, rather than for accuracy, so performing the main ``accuracy step'' in high precision and the stability steps in low precision allows to retain the stability and convergence order of the scheme while reducing the computational costs to nearly those of the low precision equivalent. They also find, however, that a naive implementation of the scheme harms the convergence order of the method and limits its accuracy.

Quantifying and also predicting the impact of using lower numerical precision is often complicated. 
For example, Prims et al.~\cite{Prims_mpc_ocean} have performed a study on two ocean models using reduced-precision emulation and a divide-and-conquer algorithm to identify -- without having to modify and recompile the code for each case -- which model variables could have their precision reduced without impacting the results. 
%
Isychev and Darulova use a domain-specific input language \cite{Isychev_2023} for testing the impact of numerical precision on a general application. Using an automatic breadth-first search, re-configuring individual functions in the executable using binary instrumentation and evaluating the impact on the results, they were able to achieve a 2x speedup on certain benchmarks without loss of accuracy.

Both methods rely on the existence of specific benchmarks or concrete scenarios, as well as of an application that can be evaluated in various precisions to determine which parts of the algorithm tolerate lower precision. 
If these parts depend on certain parameters or on the scenario, then the method must be evaluated again.
In that context, our aim is to provide a mixed-precision PDE engine for a wide set of problems, equipped with simple and easy control on which algorithmic parts should be computed in which precision. 
We also aim to provide some heuristics for how mixed precision will impact the results of these simulations, to enable users to steer their settings.

%% file: chapters/ADERV2.tex
 \section{The ADER-DG method}\label{chapter:ader}

ExaHyPE (\url{www.exahype.org}, \cite{Reinarz_ExaHyPE_2020}) solves hyperbolic PDE systems of the form 
\begin{equation}\label{eq:exa_canonical_form}
    \frac{\partial \mathbf{Q}}{\partial t} + \nabla \cdot \mathbf{F}(\mathbf{Q}, \nabla \mathbf{Q}) + \mathbf{B}(\mathbf{Q}) \cdot \nabla \mathbf{Q} = \mathbf{S}(\mathbf{Q}) + \sum_{i = 1}^{n_\text{ps}} \mathbf{\delta}_i 
\end{equation}
where $\mathbf{Q}$ is the state vector of quantities, $\mathbf{F}$ defines the conservative part and $\mathbf{B}$ defines the non-conservative part of the flux. Volume source terms are expressed in $\mathbf{S}$, and the $\mathbf{\delta}_i$ may specify $n_\text{ps}$ different point sources.
The ADER-DG method (Discontinuous Galerkin with Arbitrary high-order DERivative time stepping), as introduced by Dumbser et al.~\cite{dumbser_ader_2014}, is the main numerical approach used in ExaHyPE and also the solver used throughout this work. 
Particularly for non-linear problems, the ADER-DG solver can be complemented with an a-posteriori subcell Finite Volume limiter, whereby unphysical or oscillatory solutions are detected and solutions in respective cells are recomputed with a more robust finite volume scheme on a fine subgrid \cite{Zanotti_pos_lim}. 
In this work, however, we will restrict the discussion to scenarios that do not require limiting.
In the following, we will shortly recapitulate the key steps of the ADER-DG scheme. For this we closely follow Zanotti et al.~\cite{Zanotti_pos_lim} (see also \cite{dumbser_ader_2014,Reinarz_ExaHyPE_2020} for further details).
For simplicity, we restrict this summary to using only conservative fluxes $\mathbf{F}(\mathbf{Q})$, i.e., for solving the simplified problem $\frac{\partial \mathbf{Q}}{\partial t} + \nabla \cdot \mathbf{F}(\mathbf{Q}) = 0$. 

The spatial domain $\Omega$ is discretized on a tree-structured Cartesian grid of cells $T_i$. 
%
%
At each timestep $t^n$, the numerical solution $\mathbf{Q}_h$ is then represented within each cell $T_i$ by piecewise polynomials of maximum degree $N \geq 0$:
\begin{equation}\label{eq:piecewise_solution}
    \mathbf{Q}_h(\mathbf{x},t^n) = \sum_{l} \mathbf{\hat{Q}}_l^n \phi_l(\mathbf{x})
\end{equation}
%
with the \emph{degrees of freedom} $\mathbf{\hat{Q}}_l^n$. 
In ExaHyPE, the basis functions $\phi_l(\mathbf{x})$ are tensor-products of Lagrange polynomials of degree $N$ using Gauss-Legendre quadrature nodes as support points.
Each timestep of the ADER-DG algorithm then consists of two key steps:
\begin{enumerate}
    \item a \emph{predictor step}, which solves a space-time DG problem locally within each cell (see~section \ref{subseq:predictor});
    \item a \emph{corrector step}, which uses the results of the predictor to compute the solution at the next timestep via a discontinuous Galerkin scheme in space and considering the influence of neighbour elements (see section~\ref{subseq:corrector}).
\end{enumerate}

\subsection{Predictor}\label{subseq:predictor}

The job of the predictor is to provide for each element an estimate for the space-time solution in the next timestep, which is used in the corrector step to compute space-time approximations of the numerical fluxes at the cell boundaries (see section~\ref{subseq:corrector}). 
%
Following the usual Galerkin approach, but for a space-time element $T_i \times [t^n, t^{n+1}]$, we multiply our problem $\frac{\partial \mathbf{Q}}{\partial t} + \nabla \cdot \mathbf{F}(\mathbf{Q}) = 0$ with space-time test functions $\theta_k(\mathbf{x},t)$ and integrate in space and time, over a cell $T_i$ and a timestep $\Delta_t = [t^n, t^{n+1}]$, to obtain 
\begin{equation}\label{eq:space_time_weak_form}
    \int_{\Delta_t} \int_{T_i} \theta_k \frac{\partial \mathbf{q}_h}{\partial t} d\mathbf{x} dt + \int_{\Delta_t} \int_{T_i} \theta_k \nabla \cdot \mathbf{F}(\mathbf{q}_h) d\mathbf{x} dt = 0
\end{equation}
Here, we introduced -- analogously to \eqref{eq:piecewise_solution} -- the \emph{discrete space-time} solution $\mathbf{q}_h$ and flux $\mathbf{F}_h$: 
\begin{equation}    
    \mathbf{q}_h(\mathbf{x}, t) = \sum\limits_l \mathbf{\hat{q}_l} \theta_l(\mathbf{x}, t) 
    \quad\text{and}\quad
    \mathbf{F}_h(\mathbf{x}, t) = \sum\limits_l \mathbf{\hat{F}}_l \theta_l(\mathbf{x}, t). \label{eq:discreteFlux} 
\end{equation}
For the basis (and test) functions $\theta_l(\mathbf{x}, t)$, we again choose Lagrange polynomials with Gauss-Legendre support points. 
This choice allows the point-wise computation of the degrees of freedom for the fluxes as evaluation of the physical fluxes for the quantities $\mathbf{\hat{q}}_l$, i.e. $\mathbf{\hat{F}}_l = \mathbf{F}(\mathbf{\hat{q}}_l)$.

The time integral of the first term of \eqref{eq:space_time_weak_form} is then integrated by parts to yield:
\begin{equation}
    \int_{T_i} \theta_k(\mathbf{x}, t^{n+1}) \mathbf{q}_h d\mathbf{x} - \int_{T_i} \theta_k(\mathbf{x}, t^n) \mathbf{Q}_h d\mathbf{x} - \int_{\Delta_t} \int_{T_i} \frac{\partial \theta_k}{\partial t} \mathbf{q}_h d\mathbf{x} dt + \int_{\Delta_t} \int_{T_i} \theta_k \nabla \cdot \mathbf{F}_h d\mathbf{x} dt = 0.
\end{equation}
Substituting the discrete space-time solution and fluxes, as in equation~\eqref{eq:discreteFlux}, yields
\begin{equation} \label{eq:Picard}
    \left( \int_{T_i} \theta_k(\mathbf{x}, t^{n+1}) \theta_l(\mathbf{x}, t^{n+1}) d\mathbf{x} - \int_{\Delta_t} \int_{T_i}  \frac{\partial \theta_k}{\partial t} \theta_l d\mathbf{x} dt \right) \mathbf{\hat{q}}_l 
    + \left( \int_{\Delta_t} \int_{T_i} \theta_k \nabla \theta_l d\mathbf{x} dt \right) \mathbf{F}(\mathbf{\hat{q}}_l) 
    = \left( \int_{T_i} \theta_k(\mathbf{x}, t^n) \phi_l d\mathbf{x} \right) \mathbf{\hat{Q}}_l^n 
\end{equation}
(here using Einstein summation convention, i.e.\ summing over repeated indices, particularly $l$).
Equation~\eqref{eq:Picard} corrresponds to a nonlinear system of equations for the space-time expansion coefficients $\mathbf{\hat{q}}_l$ using the known degrees of freedom $\mathbf{\hat{Q}}_l^n$. 
In ExaHyPE this is solved in the general case via Picard iterations (essentially performing fixpoint iterations) following Dumbser et al.~\cite{Dumbser2008FiniteVS}. 
For linear PDEs, ExaHyPE provides an implementation of the Cauchy-Kowalevskaya method (see \cite{Dumbser2006BuildingBF}, e.g.) which has a much lower computational effort.

\subsection{Corrector}\label{subseq:corrector}

The corrector builds on the cell-local space-time solutions $\mathbf{q}_h(\mathbf{x}, t)$ to compute a high-order-accurate solution at timestep $t^{n+1}$, taking into account contributions from neighbour cells. 
For this, we discretise our PDE system $\frac{\partial \mathbf{Q}}{\partial t} + \nabla \cdot \mathbf{F}(\mathbf{Q}) = 0$ in space, following the Galerkin approach, and in addition integrate over one timestep $\Delta_t = [t^n, t^{n+1}]$: 
\begin{equation}\label{eq:one_step_weak_form}
    \int_{\Delta_t} \int_{T_i} \phi_k(\mathbf{x}) \frac{\partial \mathbf{Q_h}}{\partial t} d\mathbf{x} dt + \int_{\Delta_t} \int_{T_i} \phi_k(\mathbf{x}) \nabla \cdot \mathbf{F}(\mathbf{Q_h}) d\mathbf{x} dt = 0
\end{equation}
with the basis functions $\phi_k(\mathbf{x})$ of equation~\eqref{eq:piecewise_solution} as test functions. 
We then apply integration by parts on the space integral over the flux divergence term to obtain
%
\begin{equation}\label{eq:one_step_ader}
    \int_{\Delta_t} \int_{T_i} \phi_k \frac{\partial \mathbf{Q}_h}{\partial t} d\mathbf{x} dt 
    + \int_{\Delta_t} \int_{\partial T_i} \phi_k \mathbf{F}(\mathbf{Q}_h) \cdot \mathbf{n} \: dS dt 
    - \int_{\Delta_t} \int_{T_i} \nabla \phi_k \cdot \mathbf{F}(\mathbf{Q}_h) d\mathbf{x} dt = 0
\end{equation}
where $\mathbf{n}$ is the outward pointing unit normal vector on the surface $\partial T_i$ of the cell $T_i$. 
All three terms of \eqref{eq:one_step_ader} need to be integrated in time: 
for the first term, we can trivially integrate over $\frac{\partial \mathbf{Q}_h}{\partial t}$; 
for the second and third term, we perform numerical integration using the local space-time predictor $\mathbf{q}_h$ from section \ref{subseq:predictor}. 
As typical for DG schemes, the solution is discontinuous at the surface $\partial T_i$. 
For the surface integral, we therefore require a numerical flux function $\mathbf{G}$, and replace the flux term $\mathbf{F}(\mathbf{Q}_h)$ with $\mathbf{G}(\mathbf{q}_h^-, \mathbf{q}_h^+)$, using the discontinuous states on either side of $\partial T_i$.
%
We thus obtain 
\begin{equation}\label{eq_weak_form_volume_surface}
    \left( \int_{T_i} \phi_k \phi_l d\mathbf{x} \right) ( \mathbf{\hat{Q}}_l^{n+1}-\mathbf{\hat{Q}}_l^{n}) + \underbrace{\int_{\Delta_t} \int_{\partial T_i} \phi_k \mathbf{G}(\mathbf{q}_h^-, \mathbf{q}_h^+) \cdot \mathbf{n} \: dS dt}_{\textbf{surface integral}} - \underbrace{\int_{\Delta_t} \int_{T_i} \nabla \phi_k \cdot \mathbf{F}(\mathbf{q}_h) d\mathbf{x} dt}_{\textbf{volume integral}} = 0,
\end{equation}
which defines the fully discrete one-step ADER-DG scheme to compute the degrees of freedom $\mathbf{\hat{Q}}_l^{n+1}$.

The explicit scheme is stable (cf.~\cite{dumbser_2008}, e.g.) when satisfying the CFL condition
\begin{equation}\label{eq::CFL_condition}
    t^{n+1} - t^n \le \frac{C_\text{CFL} h}{d \ (2N+1) \ \lambda_\text{max}}
\end{equation}
where $C_\text{CFL}$ is a stability factor depending on the polynomial order $N$, $d$ is the number of dimensions of the problem, $h$ is the mesh size and $\lambda_\text{max}$ is the largest eigenvalue (i.e., cell-local wave speed).

Per default ExaHyPE computes the numerical fluxes $\mathbf{G}(\mathbf{q}_h^-, \mathbf{q}_h^+)$ via the so-called local-Lax-Friedrich or Rusanov flux \eqref{eq:rusanov_flux} as a Riemann solver, though other Riemann solvers can be specified by the user. With $\lambda_\text{max}^\pm$ the greatest eigenvalue of the Jacobian of the fluxes between the states on either side of the Riemann problem, the Rusanov flux is defined as:
%
\begin{equation}\label{eq:rusanov_flux}
    \mathbf{G}_\text{Rus}(\mathbf{q}^-,\mathbf{q}^+) 
    := \frac{1}{2} \left(\mathbf{F}(\mathbf{q}^-) + \mathbf{F}(\mathbf{q}^+) \right) 
    +  \frac{\lambda_\text{max}^\pm}{2} (\mathbf{q}^- - \mathbf{q}^+).
\end{equation}

%% file: chapters/Implementation.tex
\section{Implementation of mixed-precision kernels in ExaHyPE}\label{chapter:implementation}

Algorithm~\ref{alg:ader} outlines the steps of the ADER-DG method as they are implemented in ExaHyPE with key kernels marked in bold. 
\begin{algorithm}
\caption{The ADER-DG algorithm}\label{alg:ader}
\begin{algorithmic}[1]
\For {cell $T_i$ in $\Omega_h$}
    \State $ (\mathbf{q}_h, \mathbf{F}_h) \gets \text{spaceTimePredictor}(\mathbf{Q}_h)$ 
    \State $ (\mathbf{q}_h^\pm, \mathbf{F}_h^\pm) \gets \text{expansion}(\mathbf{q}_h, \mathbf{F}_h)$ \hspace{2em}
        \rlap{\smash{$\left.\begin{array}{@{}c@{}}\\{}\\{}\\{}\end{array}\color{black}\right\}%
          \color{black}\begin{tabular}{l}Predictor step, combined into the kernel \\ \textbf{fusedSpaceTimePredictorVolumeIntegral()}\end{tabular}$}}
    \State $ \mathbf{Q}_h += \text{volumeIntegral}(\partial \mathbf{F}_h))$
\EndFor
\State $\Delta t_\text{next} \gets \infty$
\For {cell $T_i$ in $\Omega$}
    \For {all faces of $T_i$}
        \State $ \mathbf{G}(\mathbf{q}_h^-, \mathbf{q}_h^+) \gets \textbf{riemannSolver}(\mathbf{q}_h^-, \mathbf{q}_h^+, \mathbf{F}_h^-, \mathbf{F}_h^+)$
            \rlap{\smash{$\left.\begin{array}{@{}c@{}}\\{}\\{}\\{}\end{array}\color{black}\right\}%
              \color{black}\begin{tabular}{l}Corrector step\end{tabular}$}}
        \State $ \mathbf{Q}_h += \textbf{faceIntegral}(\mathbf{G}(\mathbf{q}_h^-, \mathbf{q}_h^+))$
    \EndFor
    \State $\Delta t_\text{next} \gets \min(\Delta t_\text{next}, \textbf{computeTimestep}(\mathbf{Q}_h))$
\EndFor
\end{algorithmic}
\end{algorithm}
The ADER-DG algorithm consists of several operations, which can all be evaluated in various precisions: 
the space-time prediction of the cell-local polynomial, the projection of this space-time polynomial to the cell boundaries, the integration of the polynomial to compute the volume integral, the calculation of Riemann fluxes on the cell boundaries, and the integration of the Riemann fluxes to compute the cell-wise surface integral. In ExaHyPE, these operations are performed in four kernel functions:
\begin{itemize}
    \item \textbf{fusedSpaceTimePredictorVolumeIntegral()} first extrapolates the existing nodal solution into a space-time polynomial, solving equation~\eqref{eq:Picard} using either a Picard iteration or a Cauchy-Kowalevskaya procedure, then calculates the projections of this space-time polynomial to each of the cell faces. 
    Finally, it computes the contributions of the volume integral for the given cell, i.e., the term labelled \textit{volume integral} in equation \eqref{eq_weak_form_volume_surface}.
    \item \textbf{riemannSolver()} gets called once per face and computes the Riemann fluxes required by both adjacent cells, using the values on either side of the face computed by expansion(). 
    Per default this uses a Rusanov flux, but custom Riemann solvers can be implemented as well.
    \item \textbf{faceIntegral()} computes the contributions of one face to the overall surface integral of a cell. This kernel is called once per adjacent face in each cell. It corresponds to the part labelled \textit{surface integral} of equation \eqref{eq_weak_form_volume_surface}.
\end{itemize}
In each timestep, the domain is traversed twice. Once for the predictor step and a second time for the corrector step, in which the projected states and fluxes at the cell boundaries are exchanged.
The kernel \textbf{computeTimestep()} calculates the allowed timestep size for each cell. 
The uniform timestep size is computed as the smallest allowed timestep.


\subsection{ExaHyPE code generation and specifying kernel and storage precision}

ExaHyPE is built upon the parallel adaptive mesh refinement framework Peano \cite{Weinzierl:2019:Peano}. 
Peano supports flexible implementation of data structures through the tool DaStGen~\cite{bungartz_dastgen}, which allows the generation of C++ classes from a python API. 
These classes can be expanded with additional attributes, for which getter and setter functions are automatically created. 
DaStGen also adds utility functions, such as support for MPI datatypes as well as data exchanges between adjacent objects. 
The latest version of the tool, DaSTGen2, also supports optimized storage via non-IEEE floating-point format data. 
The generated C++ code can be expanded with compiler annotations that specify how to efficiently compress and store the data. 
In current work, this is extended towards memory-centric extensions to C++ \cite{PawelMemoryAccessor}.
Our mixed-precision implementation heavily relies on this flexibility of data structures. 

For the implementation of all ADER-DG kernels, ExaHyPE uses a role-oriented code generation approach~\cite{gallard:2020}. 
All kernel functions are generated from Jinja2~\cite{jinja2_docs} templates, taking into account the specifics of the problem to be solved. 
For example the code can expose exact memory requirements, loop counts or array lengths (depending on the number of quantities or polynomial degree, e.g.) to the compiler, which allows for more efficient optimization (esp.\ vectorization). 
In addition the need for branching is minimized: while static code might need to check on every execution whether or not to evaluate certain kernels, the code generator eliminates unnecessary calls upon generation.

To specify the precision to store persistent and temporary data, as well as to compute the various kernels of ADER-DG, we provide the following four control parameters to the code generation framework: 
\begin{description}
\item[predictor] -- the precision in which the space-time predictor step, i.e., the kernel \textbf{fusedSpaceTimePredictorVolumeIntegral()}, will be computed and stored, 
\item[corrector] -- the precision in which the information on the faces will be stored and the corrector step, composed of the kernels \textbf{riemannSolver()} and \textbf{faceIntegral()}, will be computed,
\item[storage] -- the precision in which the persistent solution will be stored between grid traversals, and
\item[Picard] -- the precision in which the Picard iterations within the space-time preditor are computed 
   (only for nonlinear PDE systems; for linear PDE systems, the space-time prediction is computed via the CK procedure which is directly integrated into the predictor step). 
\end{description}
Users can set any of these modes to any of the floating point data types offered by C++: in our case \texttt{bfloat16}, \texttt{fp16}, \texttt{float}, \texttt{double} and \texttt{long double} (compare table \ref{table:precision_bits}).
If not specified, the modes default to using \texttt{double} precision.


\subsection{Implementing mixed-precision kernels in ExaHyPE's code generation framework}

ExaHyPE2 combines pre-compiled static kernels with application-specific (i.e., PDE-specific) code. The latter can be further differentiated between user functions (implemented for a particular PDE system) and generated code.
All of these may be required in one or more precisions, but for each of these different design principles take precedence:
\begin{itemize}
    \item Static kernels are only written once and reused between all applications. They are implemented with minimal information about the final application, which they must be able to serve regardless of the final precision. The main goals for these kernels are maximising efficiency while keeping the code maintainable.
    \item Generated kernels are instantiated when the actual application is defined -- from pre-existing Jinja2 templates and only in the required variants. This allows for further optimisation, but means that not all variants can be tested in advance. To simplify software maintenance, these variants should remain consistent with one another and code duplication should be minimized. 
    For example, a single Jinjia2 template should allow for instantiation into any given version, instead of having separate templates for any possible combination of parameters.
    \item User functions are the single point of contact of an application developer (``user'') with the code. 
    Examples include functions to compute initial and boundary conditions of the problem, as well as eigenvalues or fluxes of the underlying PDE system.
    As such, user functions must be easily understandable and the effort for the user should be minimized.
    Users who do not request mixed-precision should only be presented with a single precision. 
    Users who request multiple precision should not have to implement their functions multiple times.
\end{itemize}

As the target language for all generated code is C++, one option to provide kernels in different precision would be to use C++ templates and let the compiler instantiate one specialisation for each variant. 
In that case, the compiler needs to be aware of each required specialisation at compile time.
This would, however, conflict with ExaHyPE's design as an engine: while generated kernels are generated on demand from user specifications of their application, static kernels (such as mesh traversal, I/O, \dots) are application-independent. 
The respective source code is pre-compiled into object files in order to speed up the compilation of individual applications. 
To still provide any possible specialisation a user might request, every possible specialisation of every function would have to be pre-compiled. This would strongly increase initial compilation times and object file sizes of the project.

The standard way of resolving this, which is used in ExaHyPE for all static kernels, is by providing template definitions through header files that are then included in every source file which makes use of them. This means that during compilation, the corresponding translation unit can generate whichever template specialisation it requires. This however increases both the time required for compilation and the size of the resulting object files since the same template specialization will be compiled multiple times and stored within each object file that requires it. 
For the sake of readability, template definitions in ExaHyPE are implemented in separate \textit{.cpph} files, which are included by their corresponding header file. This allows for clearer separation of definitions and declarations.

\begin{figure}
        \begin{lstlisting}[language=Python,frame=lrtb]
{% for precision_num in range(0,computePrecisions|length) %}
template void {{codeNamespace}}::faceIntegral(
  {{computePrecisions[precision_num]}} *lduh, 
  const {{correctorComputePrecision}} *const lFhbnd,
  const int direction, const int orientation,
  const double inverseDxDirection
);
{% endfor %}
    \end{lstlisting}
    \caption{Example how to generate C++ template specialisations in all required precisions using the Jinja2 template engine.}
\end{figure}

In contrast to the static kernels, application-specific generated kernels are only generated once the user has fully defined the problem, such that the precisions required at runtime are known. 
Therefore, explicit C++ template specialisations can be employed instead of providing template definitions through header files.
This has the main benefit that if any changes need to be made to other sections of the application without needing to regenerate the kernels, such as changing the boundary conditions of the problem, one would not need to re-compile the templated functions. It also makes the code more legible and maintains a more thorough separation between the sections of the code defining the algorithm and those which only make use of it.

Finally, all user functions are provided only in the precisions explicitly requested for the actual application. 
If the user specifies to perform all computations in uniform precision, the user functions don't need to be templates and are therefore provided only once, in fixed precision. 
Otherwise, the function stubs are generated as templates (see figure~\ref{fig:flux_stub}) together with any required template specialisations.

\begin{figure}
    \begin{minipage}[t]{.48\textwidth}
        \begin{lstlisting}[language=C++,frame=lrtb]

void flux(
  const double* __restrict__ Q,
  double                     t,
  double                     dt,
  int                        normal,
  double* __restrict__       F
) override;

        \end{lstlisting}%
    \end{minipage}%
    \hfill
    \begin{minipage}[t]{.48\textwidth}
        \begin{lstlisting}[language=C++,frame=lrtb]
template<typename T>
void flux(
  const T* __restrict__ Q,
  double                t,
  double                dt,
  int                   normal,
  T* __restrict__       F
);
        \end{lstlisting}%
    \end{minipage}%

    \caption{Different declaration of the user function \texttt{flux()}, depending on whether the user has specified multiple precisions (right) or to perform all computations in uniform \texttt{double} precision (left).
    Note that \texttt{flux()} is called by the predictor and the corrector.}
    \label{fig:flux_stub}
\end{figure}

The precision in which persistent data is stored is also controlled via Jinja2 templates. 
This data is stored as values inside of a \texttt{struct}, which allows for some associated functions for operations such as initializing and deleting the data, or merging data from two neighbouring faces. This makes the usage of different precisions particularly simple since the architecture will only be using the containing \texttt{struct} in most situations.
Only operations that interact with the data contained within this \texttt{struct} need take that particular data into account.

Finally, the ADER-DG kernels rely on several precomputed matrices for many of their operations. 
These matrices contain information such as the position of the nodes of the Lagrange-polynomial basis functions, the weights for Gaussian quadrature, the reference element stiffness matrix or the weights to compute the projection of the solution to the faces of a cell.
For simplicity, one might be tempted to only store these in the highest precision. 
However, C++ always performs operations in the highest precision among the data used: 
Multiplying a \texttt{float} with a \texttt{double} and storing the result as a \texttt{float} would cast the \texttt{float} to \texttt{double}, perform a \texttt{double}-precision multiplication and cast the result to \texttt{float}. 
This would not perform the intended \texttt{float}-multiplication, and thus not profit from a longer vectorisation length, for example. 
We therefore need to provide all precomputed matrices in every required precision. 
The matrices are stored in two classes, one for data related to Gaussian quadrature, such as quadrature weights or nodes, and one for data associated with ADER-DG, such as the stiffness or projection matrices. 
In that way, these constant matrices are shared between multiple instances of the functions, but initialized only once. 
Analogously to the kernels, the classes are templated in every required precision.
We can guarantee that the kernels have access to every required version of the matrices simply by initializing these once in every precision required by the kernels.

%% file: chapters/Scenarios.tex
\section{Test scenarios}\label{chapter:verification}

In this chapter, we outline five scenarios with known analytical solution to evaluate the accuracy of results when using low- or mixed-precision ADER-DG. 
The scenarios are based on four common hyperbolic PDE systems: the acoustic and elastic wave equations, the Euler equations and the shallow water equations.
We evaluate the impact of rounding errors on convergence properties and on specific numerical challenges, and compare these to discretisation errors. 
For the non-linear scenarios, we have chosen scenarios that do not require a limiter for the ADER-DG method. 

\subsection{Linear acoustic wave equation} \label{subseq:lin_acoustic}

The linear acoustic wave equation (cmp., e.g., LeVeque \cite{leveque_2002}) describes the propagation of acoustic waves in a homogeneous medium, for example in a fluid under the assumption that viscosity and heat conduction are negligible. 
Formulated via the velocity $v$ and the pressure $p$ of the fluid, we obtain the following PDE system (in 2D):

\begin{equation} \label{eq:acoustic}
    \begin{aligned}
        \frac{\partial p}{\partial t} + K \left( \frac{\partial v_x}{\partial x} + \frac{\partial v_y}{\partial y} \right)&= 0 \\
        \rho \frac{\partial v_x}{\partial t} + \frac{\partial p}{\partial x} &= 0 \\
        \rho \frac{\partial v_y}{\partial t} + \frac{\partial p}{\partial y} &= 0
    \end{aligned}
\end{equation}
(with $\rho$ the density of the fluid, $K$ its bulk modulus).
To match it to the canonical form \eqref{eq:exa_canonical_form}, we write the PDE system as
\begin{equation}\label{eq:lin_acoustic}
     \frac{\partial \mathbf{Q}}{\partial t} + \nabla \cdot \mathbf{F}(\mathbf{Q}, \nabla \mathbf{Q}) = 
     \frac{\partial \mathbf{Q}}{\partial t} + \mathbf{A} \frac{\partial \mathbf{Q}}{\partial x} + \mathbf{B} \frac{\partial \mathbf{Q}}{\partial y} 
     = 0
\end{equation}
with 
$\mathbf{Q} =
    \begin{pmatrix} 
        p & v_x & v_y
    \end{pmatrix}^T
$ the vector of quantities and matrices 
\begin{equation}\label{eq:acoustic_matrices}
    \mathbf{A} =
    \begin{pmatrix} 
        0 & K & 0 \\
        \frac{1}{\rho} & 0 & 0 \\
        0 & 0 & 0
    \end{pmatrix}
    \qquad
    \mathbf{B} =
    \begin{pmatrix} 
        0 & 0 & K \\
        0 & 0 & 0 \\
        \frac{1}{\rho} & 0 & 0
    \end{pmatrix}.
\end{equation}
%
As verification scenario for the acoustic wave equation, we compute a planar wave propagating with the \textbf{acoustic wave speed} $c = \sqrt{\frac{K}{\rho}}$ (see section \ref{eq:elastic} for details).

\subsection{Linear elastic wave equations}  \label{eq:elastic}

The linear elastic wave equations \cite{leveque_2002} describe the propagation of waves through heterogeneous media assuming a linear relationship between stress and strain. Formulated via the velocity $v$ and the stress $\sigma$, we obtain the following PDE system:
\begin{equation}\label{eq:equations_of_motion}
    \begin{aligned}
        \frac{\partial \sigma_{xx}}{\partial t} = (\lambda + 2\mu) \frac{\partial v_x}{\partial x} + \lambda \frac{\partial v_y}{\partial y} \\
        \frac{\partial \sigma_{yy}}{\partial t} = \lambda \frac{\partial v_x}{\partial x} + (\lambda + 2\mu) \frac{\partial v_y}{\partial y} \\
        \frac{\partial \sigma_{xy}}{\partial t} = \mu \left(\frac{\partial v_x}{\partial y} + \frac{\partial v_y}{\partial x}\right) \\
        \rho \frac{\partial v_x}{\partial t} = \frac{\partial  \sigma_{xx}}{\partial x} + \frac{\partial \sigma_{xy}}{\partial y} \\
        \rho \frac{\partial v_y}{\partial t} = \frac{\partial  \sigma_{xy}}{\partial x} + \frac{\partial \sigma_{yy}}{\partial y}.
    \end{aligned}
\end{equation}
%
Similarly to the acoustic wave equation \eqref{eq:lin_acoustic}, this can be rewritten in the form
\begin{equation}\label{eq:lin_elastic}
    \frac{\partial \mathbf{Q}}{\partial t} + \nabla \cdot \mathbf{F}(\mathbf{Q}, \nabla \mathbf{Q}) = 
    \frac{\partial \mathbf{Q}}{\partial t} + 
    \mathbf{A} \frac{\partial \mathbf{Q}}{\partial x} + \mathbf{B} \frac{\partial \mathbf{Q}}{\partial y}
    = 0
\end{equation}
with 
$\mathbf{Q} =
    \begin{pmatrix} 
        \sigma_{xx} & \sigma_{yy} & \sigma_{xy} & v_x & v_y
    \end{pmatrix}^T
$ and matrices 
{\small
\begin{equation}\label{eq:elastic_matrices}
    \mathbf{A} =
    \begin{pmatrix} 
        0 & 0 & 0 & - ( \lambda + 2\mu ) & 0 \\
        0 & 0 & 0 & - \lambda & 0 \\
        0 & 0 & 0 & 0 & -\mu \\
        - \frac{1}{\rho} & 0 & 0 & 0 & 0 \\
        0 & 0 & - \frac{1}{\rho} & 0 & 0
    \end{pmatrix} \qquad 
    \mathbf{B} =
    \begin{pmatrix} 
        0 & 0 & 0 & 0 & - \lambda \\
        0 & 0 & 0 & 0 & - ( \lambda + 2\mu ) \\
        0 & 0 & 0 & -\mu & 0 \\
        0 & 0 & - \frac{1}{\rho} & 0 & 0 \\
        0 & - \frac{1}{\rho} & 0 & 0 & 0
    \end{pmatrix}. 
\end{equation}
}
$\lambda$ and $\mu$ are the so-called Lamé parameters, which describe the material properties of the medium. In particular they determine the \textbf{elastic wave speeds} $c_p = \sqrt{\frac{\lambda + 2\mu}{\rho}}$ and $c_s = \sqrt{\frac{\mu}{\rho}}$ for pressure waves and shear waves, respectively.

\subsubsection{Planar waves}\label{subsubseq:planar_waves}

As verification scenario for both the acoustic and elastic wave equation, we compute a planar wave moving through the domain. These are frequently used for convergence studies of high order methods \cite[e.g.]{dumbser_ader_2006,wolf_aderdg_seismic}.
%
The planar wave scenario uses a sine mode propagated in direction $\mathbf{k} = (k_x, k_y)^T$ with frequency $\omega$ and a state vector $\mathbf{q}_0$ of amplitudes as initial condition, where $\omega$ and $\mathbf{q}_0$ are an eigenpair of the matrix $k_x \mathbf{A} + k_y \mathbf{B}$, with $\mathbf{A}$ and $\mathbf{B}$ the matrices from \eqref{eq:acoustic_matrices} or \eqref{eq:elastic_matrices}. This leads to following analytical solution:
\begin{equation}
    \mathbf{Q}(\mathbf{x},t) =
    \omega \sin( \omega t - \mathbf{k} \cdot \mathbf{x} ) \mathbf{q}_0
\end{equation}

For simplicity we choose the domain $[-1,1]\times [-1,1]$ and apply periodic boundary conditions. For the acoustic equation we set the bulk modulus $K=4$ and density $\rho=1$, which leads to a wave speed of $c=2$, and a wave travelling in diagonal direction, such that the solution returns to its initial condition at every $t = i\sqrt{2}, i= 1, 2, \dots$. 
For the elastic wave equation we choose $\lambda=2$, $\mu=1$ and $\rho=1$, which leads to wave speeds of $c_p=2$ and $c_s=1$, and a superposition of two waves, each travelling in a cardinal direction, such that the solution returns to its initial condition at every $t = i$ for $i= 1, 2, \dots$. 
We simulate two domain traversals of the wave for each equation.

\subsection{Euler equations}

The Euler equations describe the flow of a compressible fluid given the assumption that both heat conduction and viscosity are negligible. Formulated via the density $\rho$, the momentum $\rho v$ and the total energy $E$ they can be described by the following PDE system \cite{dumbser_2008}:
\begin{equation}
    \begin{aligned}
    \frac{\partial \rho}{\partial t} + \frac{\partial \: \rho v_x}{\partial x} + \frac{\partial \: \rho v_y}{\partial y}  &= 0 \\
    \frac{\partial \: \rho v_x}{\partial t} + \frac{\partial \: \rho v_x^2}{\partial x} + \frac{\partial \: \rho v_x v_y}{\partial y} + \frac{\partial p}{\partial x} &= 0 \\
    \frac{\partial \: \rho v_y}{\partial t} + \frac{\partial \: \rho v_x v_y}{\partial x} + \frac{\partial \: \rho v_y^2}{\partial y} + \frac{\partial p}{\partial y} &= 0 \\
    \frac{\partial E}{\partial t} + \frac{\partial \: v_x(\rho E + p)}{\partial x} + \frac{\partial \: v_y (\rho E + p)}{\partial y} &= 0 
    \end{aligned}
\end{equation}
which can be written in the form
\begin{equation}\label{eq:euler_eq}
    \frac{\partial \mathbf{Q}}{\partial t} + \mathbf{F}(\mathbf{Q}, \nabla \mathbf{Q}) =
    \frac{\partial \mathbf{Q}}{\partial t} + 
    \frac{\partial \mathbf{G}(\mathbf{Q})}{\partial x} + \frac{\partial \mathbf{H}(\mathbf{Q})}{\partial y}
    = 0
\end{equation}
with 
$\mathbf{Q} =
    \begin{pmatrix} 
        \rho & \rho v_x & \rho v_y & E
    \end{pmatrix}^T
$ the vector of conserved quantities, and with flux vectors 
\begin{equation}
    \mathbf{G}(\mathbf{Q}) =
    \begin{pmatrix}
        \rho v_x & \rho v_x^2 + p & \rho v_x v_y & v_x (\rho E+p) \:
    \end{pmatrix}^T 
    \quad\text{and}\quad
    \mathbf{H}(\mathbf{Q}) =
    \begin{pmatrix}
        \rho v_y & \rho v_x v_y & \rho v_y^2 + p & v_y (E+p) \:
    \end{pmatrix}^T. 
\end{equation}
For an ideal gas, the pressure $p$ is related to energy $E$ and velocity $v$ via the equation of state 
\begin{equation}
    p = (\gamma -1) (E - \textstyle\frac{1}{2} \rho v^2 )
\end{equation}
with $\gamma$ the ratio of specific heats.
The characteristic wave speeds of the Euler equations are $v+\sqrt{\frac{\gamma p}{\rho}}$ and $v-\sqrt{\frac{\gamma p}{\rho}}$.

\subsubsection{Advection of a Gaussian density bell}\label{subsubseq:gaussian_bell}

As the first verification scenario for the Euler equations, we use the \textit{advection of a smooth density bell in 2D}, as used by Ioriatti et al.~\cite{dumbser_euler} for evaluating high-order discontinuous Galerkin schemes.
The test prescribes a Gaussian initial density distribution (cmp.~figure \ref{fig:gaussian_bell})
\begin{align*}
    \rho(x,y) = \rho_0 \bigl( 1 + e^{-50\sqrt{x^2+y^2}} \bigr)
\end{align*}
that is advected along with a fluid of constant pressure $p=1$ and velocity $v_x = v_y = 1$. 
With $\rho_0 = 0.02$, the initial conditions become:
\begin{equation}
    \mathbf{Q}(x,y) =
    \begin{pmatrix} 
        0.02 ( 1 + e^{-50\sqrt{x^2+y^2}} ) \\ 
        0.02 ( 1 + e^{-50\sqrt{x^2+y^2}} ) \\ 
        0.02 ( 1 + e^{-50\sqrt{x^2+y^2}} ) \\ 
        \frac{1}{\gamma-1} + 0.02 ( 1 + e^{-50\sqrt{x^2+y^2}} )
    \end{pmatrix}
\end{equation}
The computational domain is $[-1,1]\times[-1,1]$, periodic boundary conditions are used and the total simulated time is $t_\text{end} = 2$, or two full traversals of the domain by the density bell.
Additionally we choose $\gamma=1.4$, which is the ratio of specific heats for dry air at room temperature.

\begin{figure}
    \begin{minipage}[b]{0.48\linewidth}
        \input{figures/analytical_scenarios/gaussian_bell}
        \caption{Initial density $\rho(x,y) = 0.02 \left(1+e^{-50(x^2+y^2)}\right)$ of the Gaussian bell scenario for the Euler equations.}
        \label{fig:gaussian_bell}
    \end{minipage}
    \hfill
    \begin{minipage}[b]{0.48\linewidth}
        \input{figures/analytical_scenarios/isentropic_vortex}
        \caption{Initial density $\rho(x,y) = \rho_{\infty} + \delta \rho$ of the static isentropic vortex scenario for the Euler equations.}
        \label{fig:isentropic_vortex}
    \end{minipage}%
\end{figure}

\subsubsection{Static isentropic density vortex}\label{subsubseq:isentropic_vortex}

As a second scenario, we use the isentropic density vortex, as described by Hu and Shu \cite{Shu_vortex}, prescribed by the initial conditions
\begin{equation}
    \mathbf{Q}(x,y) =
    \begin{pmatrix} 
        \rho_{\infty} + \delta \rho \\ 
        ( \rho_{\infty} + \delta \rho ) ( u_{\infty} + \delta u ) \\ 
        ( \rho_{\infty} + \delta \rho ) ( u_{\infty} + \delta v ) \\ 
        \frac{p_{\infty} + \delta p }{\gamma-1} + 0.5( \rho_{\infty} + \delta \rho )(( u_{\infty} + \delta u )^2 + ( u_{\infty} + \delta v )^2)
    \end{pmatrix}
    \quad\text{with}\quad
    \begin{dcases}
        \delta \rho  = ( 1 + \delta T )^{\frac{1}{\gamma - 1}} - 1 \\
        \delta p = ( 1 + \delta T )^{\frac{1}{\gamma - 1}} - 1 \\
        \delta T = \frac{ ( 1 - \gamma ) \beta^2 }{8 \gamma \pi^2} e^{1-(x^2+y^2)} \\
        \delta u = -y \frac{\beta}{2 \pi} e^\frac{1-(x^2+y^2)}{2} \\
        \delta v = x \frac{\beta}{2 \pi} e^\frac{1-(x^2+y^2)}{2}
    \end{dcases}
\end{equation}
The solution is a vortex that does not deform, but simply shifts through the domain at constant velocity $[u_{\infty}, v_{\infty}]$. As the name implies, the total entropy $S=\frac{p}{\rho^\gamma}=1.0$ remains constant in time and over the entire domain.
We set as parameters $\rho_{\infty}=1.0$, $p_{\infty}=1.0$ and $\beta=5.0$, and as velocities $u_{\infty} = v_{\infty} = 0$ to obtain a static solution over time: 
at each point of the vortex, the balance of all fluxes is $0$.  
In practice, however, the isentropic vortex scenario is unstable, since any slight error causes a deformation of the vortex. 
These errors grow, as resultant fluxes no longer compensate, which in turn causes the vortex deformation to grow. 
As such the isentropic vortex is a very sensitive test for the effects of numerical precision, since the analytic solution is static, but small numerical errors accumulate and  cause the simulations to diverge from its initial conditions.
For our tests, we use the computational domain $[-5,5]\times[-5,5]$, 
apply periodic boundary conditions and simulate until $t_\text{end}=10.0$.

\subsection{Shallow water equations} \label{sec:swe}

The shallow water equations describe the behavior of fluids in situations where quantities may be averaged in the vertical direction. This is the case if the depth of the fluid is very small (when compared to its horizontal extent) -- hence the name \textit{shallow} water equations -- or more general if the characteristic wave lengths are large compared to the water depth. 
Neglecting viscosity and any other source terms (Coriolis forces, e.g.), we formulate the shallow water equations via the water height $h$, the discharge $h v$ and the bathymetry $b$ (compare \cite[e.g.]{LeVeque_George_Berger_2011}): 
\begin{equation} \label{eq:swe}
    \begin{aligned}
        \frac{\partial h}{\partial t} + \frac{\partial h v_x}{\partial x} + \frac{\partial h v_y}{\partial y} &= 0 \\
        \frac{\partial h v_x}{\partial t} = \frac{\partial h v_x^2}{\partial x} + \frac{\partial \: h v_x v_y}{\partial y} + h g \frac{\partial (b + h)}{\partial x} &= 0 \\
        \frac{\partial h v_y}{\partial t} = \frac{\partial \: h v_x v_y}{\partial x} + \frac{\partial h v_y^2}{\partial y} + h g \frac{\partial (b + h)}{\partial y} &= 0
    \end{aligned}
\end{equation}
Again, to match it to the canonical form \eqref{eq:exa_canonical_form}, we write the PDE system as
%
\begin{equation}\label{eq:swe_canonical}
    \frac{\partial \mathbf{Q}}{\partial t} + \nabla \cdot \mathbf{F}(\mathbf{Q}, \nabla \mathbf{Q}) + B(\mathbf{Q}) \cdot \nabla \mathbf{Q} = \frac{\partial \mathbf{Q}}{\partial t} + \frac{\partial \mathbf{G}(\mathbf{Q})}{\partial x} + \frac{\partial \mathbf{H}(\mathbf{Q})}{\partial y} + \mathbf{C(Q)} \frac{\partial \mathbf{Q}}{\partial x} + \mathbf{D(Q)} \frac{\partial \mathbf{Q}}{\partial y} = 0
\end{equation}
with 
$\mathbf{Q} =
    \begin{pmatrix} 
        h & h v_x & h v_y & b
    \end{pmatrix}^T
$ the vector of quantities, conservative fluxes $F(\mathbf{Q}, \nabla \mathbf{Q})$ defined via 
\begin{equation} \label{eq:swe_cons}
    \mathbf{G}(\mathbf{Q}) =
    \begin{pmatrix}
        v_x & h v_x^2 & h v_x v_y & 0 \:
    \end{pmatrix}^T
    \quad\text{and}\quad
    \mathbf{H}(\mathbf{Q}) =
    \begin{pmatrix}
        v_y & h v_x v_y & h v_y^2 & 0 \:
    \end{pmatrix}^T \\
\end{equation}
and non-conservative fluxes $B(\mathbf{Q}) \cdot \nabla \mathbf{Q}$ resulting from matrices 
\begin{equation} \label{eq:swe_noncons}
    \mathbf{C(Q)} =
    \begin{pmatrix} 
        0 & 0 & 0 & 0 \\
        gh & 0 & 0 & gh \\
        0 & 0 & 0 & 0 \\
        0 & 0 & 0 & 0
    \end{pmatrix}
    \quad\text{and}\quad
    \mathbf{D(Q)} =
    \begin{pmatrix} 
        0 & 0 & 0 & 0 \\
        0 & 0 & 0 & 0 \\
        gh & 0 & 0 & gh \\
        0 & 0 & 0 & 0
    \end{pmatrix}.
\end{equation}
Note that this formulation adds the equation $\frac{\partial b}{\partial t} = 0$ to the PDE system \eqref{eq:swe}. 
Thus, the bathymetry $b$ is formally treated as a variable quantity, which however stays constant over the entire simulation. 

We consider the shallow water equations in the context of tsunami simulation \cite[e,g,]{LeVeque_George_Berger_2011}.
The characteristic wave speeds of the shallow water equations, given by $v+\sqrt{gh}$ and $v-\sqrt{gh}$ \cite[e.g.]{leveque_2002}, then determine the propagation speed of the tsunami.

\subsubsection{Lake at rest}\label{subsubseq:resting_lake}

The lake at rest scenario prescribes zero initial velocities and constant water surface elevation $\eta := h+b \stackrel{!}{=} \eta_0$ within the entire domain, but with variable bathymetry $b(x,y)$. 
We use a sinusoidal bathymetry, such that the initial conditions are:
\begin{equation} \label{eq:lake_initial}
    \begin{aligned}
        b (x,y) &= \sin(2\pi (x+y)) \\
        h (x,y) &= \eta_0 - b(x,y) \\
        v_x(x,y) &= v_y(x,y) = 0
    \end{aligned}
\end{equation}
It is easy to verify that the initial conditions \eqref{eq:lake_initial} provide a constant-in-time solution (a ``lake at rest'') to the PDE system~\eqref{eq:swe}.  
Numerical schemes need to reflect that property as accurately as possible, i.e., need to be \emph{well-balanced}: 
For tsunami simulations, e.g., the ocean depth (i.e., $h$ but also $-b$) is orders of magnitude larger than the sea surface elevation $\eta = h+b$, such that even small errors in the treatment of the $b$- and $h$-terms in equation~\eqref{eq:swe} can lead to errors that are larger than the elevation $\eta$.
Well-balanced schemes must preserve the constant surface $\eta_0$ exactly \cite[e.g.]{dumbser_hllem,Gosse_well_balancedness}.
Only rounding errors may occur when evaluating the scheme with finite precision. 
The lake at rest scenario is thus an important benchmark for the numerical accuracy of shallow water solvers. 

Discontinuous Galerkin schemes need to satisfy the well-balanced property for the volume integrals (cf.\ Eq.~\eqref{eq_weak_form_volume_surface}) as well as for the numerical fluxes between the discontinuous cell-local approximations of two elements. 
Our formulation \eqref{eq:swe_canonical} of the shallow-water equations leads to a well-balanced discretisation within each cell, as the conservative fluxes \eqref{eq:swe_cons} as well as the non-conservative fluxes \eqref{eq:swe_noncons} become zero for constant water level $\eta=h+b$ and zero velocities $v_x=v_y=0$. 
To ensure that numerical fluxes also lead to a well-balanced scheme, we use a modified version of the Rusanov flux, analogous to \cite{Leo_thesis}:
\begin{equation}\label{eq:swe_custom_riemann}
    \mathbf{F}_\text{Rusanov} = \frac{1}{2} (\mathbf{F}^- + \mathbf{F}^+ ) - \frac{1}{2} \mathbf{B}(\mathbf{q}^-, \mathbf{q}^+) \cdot \mathbf{n} + \frac{1}{2} \begin{pmatrix}
        h^- + b^- - h^+ - b^+ \\ v_{x}^- - v_{x}^+ \\ v_{y}^- - v_{y}^+ \\ 0 \:
    \end{pmatrix}
\end{equation}
where $\mathbf{F}^-$ and $\mathbf{F}^+$ are the fluxes on either side of the discontinuity, $\mathbf{B}(\mathbf{q}^-, \mathbf{q}^+)$ is the non-conservative flux between the states on either side, $\mathbf{q}^-$ and $\mathbf{q}^+$, and $\mathbf{n}$ is the direction relative to the cell face. This is a similar principle to that of \cite{dumbser_hllem}, using an intermediate state between the left and right states to evolve the wave.

The constant surface elevation $\eta_0$ in the initial condition~\eqref{eq:lake_initial} is typically set to $0$.
However, we purposefully use a non-zero value for the elevation, for example $\eta_0 = 2$, as this has a noticeable effect on rounding errors. 
While this shift of water level has no impact on the well-balancedness of the numerical scheme, we obtained much better numerical results when setting $\eta_0 = 0$. 
In fact, numerical experiments for $\eta_0 = 0$ caused no artificial waves, regardless of the floating-point precision. 
This is different for $\eta_0 \not= 0$ (see section~\ref{subseq:static}). 
And note that even for simulations that follow the convention of $\eta_0 = 0$, well-balancedness for non-zero water elevation is relevant, as it applies whenever there is a smooth surface elevation that extends over non-constant bathymetry $b$ (consider a tsunami with long wave length, e.g.).

\subsection{On the computation of errors}\label{subseq:error_estimate}  

For many of the evaluation scenarios, we know the analytical solution $\mathbf{Q}(\mathbf{x},t)$, which suggests computing the error in the numerical solution $\mathbf{Q}_h(\mathbf{x},t^n)$ at a timestep $t^n$ via a suitable norm -- such as using the $L^2$ norm, 
\begin{equation} \label{eq:L2exact}
       \| \mathbf{Q} - \mathbf{Q}_h \|_{2} = \sqrt{\int_\Omega \sum_v \bigl( \mathbf{Q}_v(\mathbf{x},t^n) - \mathbf{Q}_{h,v}(\mathbf{x},t^n) \bigr)^2 \,d\mathbf{x}},
\end{equation}
or some other suitable norm (maximum norm, etc.). 
In practice, the computation of \eqref{eq:L2exact} requires evaluation of $\mathbf{Q}$ and $\mathbf{Q}_h$ at integration points $\mathbf{x}_i$, leading to some integration rule (with integration weights $w_i$) 
\begin{equation} \label{eq:L2integral}
       \| \mathbf{Q} - \mathbf{Q}_h \|_{2} \approx \sqrt{ \sum_{c \in \text{cells}} \sum_{i \in \text{nodes}} \sum_v V_c w_i (\mathbf{Q}_v(\mathbf{x}_i,t^n) - \mathbf{Q}_{h,v}(\mathbf{x}_i,t^n))^2}
\end{equation}
summing over all DG elements $c$ (with volume $V_c$). 
As indicated in equation \eqref{eq:piecewise_solution}, our DG discretization uses tensor products of Lagrange polynomials as nodal basis functions. These use Gaussian quadrature nodes $\mathbf{x}_i$ as support.
The coefficients of the element-wise DG polynomial $\mathbf{Q}_h^\ast(\mathbf{x},t^n)$ that interpolates the exact solution at the quadrature nodes, are thus given as $\mathbf{Q}_h^\ast(\mathbf{x}_i,t^n)$. 
And due to Gaussian quadrature (with respective weights $w_i$) we can efficiently and with good accuracy compute the norm
\begin{equation} \label{eq:L2gauss}
       \| \mathbf{Q}_h^\ast - \mathbf{Q}_h \|_{2} = \sqrt{ \sum_{c \in \text{cells}} \sum_{i \in \text{nodes}} \sum_v V_c w_i (\mathbf{Q}_{h,v}^\ast(\mathbf{x}_i,t^n) - \mathbf{Q}_{h,v}(\mathbf{x}_i,t^n))^2}.
\end{equation}
For comparison, to determine the maximum difference 
$\bigl| \mathbf{Q}(\mathbf{x},t^n) - \mathbf{Q}_h(\mathbf{x},t^n) \bigr|$, 
or alternatively
$\bigl| \mathbf{Q}_h^\ast(\mathbf{x},t^n) - \mathbf{Q}_h(\mathbf{x},t^n) \bigr|$, 
in each element, it would not suffice to evaluate $\mathbf{Q}$ (or $\mathbf{Q}_h^\ast$) and $\mathbf{Q}_h$ at integration points, but subsampling on a finer grid and evaluation of $\mathbf{Q}$ (or $\mathbf{Q}_h^\ast$) and $\mathbf{Q}_h$ at element boundaries would be necessary. 
Also, the maximum error is more prone to masking effects: in the case of shocks or steep gradients, for example, such local errors might dominate the error more strongly than for the $L^2$ norm.

We therefore take $\mathbf{Q}_h^\ast(\mathbf{x},t^n)$ as the \textit{reference solution} and, unless otherwise noted, compute the \textit{L2 error} $\| \mathbf{Q}_h^\ast - \mathbf{Q}_h \|_{2} $, evaluated following \eqref{eq:L2gauss}, to quantify the error in the numerical solution $\mathbf{Q}_h$. 
This implies that $\mathbf{Q}_h$ may have been computed using lower or mixed floating-point precision and of course also using different mesh size or polynomial order. 
We always evaluate \eqref{eq:L2gauss} in double precision, i.e, we cast any low-precision coefficients of $\mathbf{Q}_h$ to \texttt{fp64} before evaluation.
We thus ignore any error $\| \mathbf{Q} - \mathbf{Q}_h^\ast \|_{2}$ that results from not being able to represent the exact solution $\mathbf{Q}(\mathbf{x},t^n)$ via element-wise DG-polynomials -- including that we use $\mathbf{Q}_h^\ast(\mathbf{x},0)$ instead of $\mathbf{Q}(\mathbf{x},0)$ as starting condition for a time-dependent DG-solution. 
However, this contribution to the overall errors is independent of the applied floating-point precision.

%% file: figures/analytical_scenarios/gaussian_bell.tex
\begin{tikzpicture}{}
    \begin{axis}
        \addplot3[surf, domain=-1:1, samples=70, samples y=70] function {0.02*(1+exp(-0.5*(x*x+y*y)/0.01))};
    \end{axis}
\end{tikzpicture}

%% file: figures/analytical_scenarios/isentropic_vortex.tex
\begin{tikzpicture}{}
    \begin{axis}
        \addplot3[surf, domain=-3:3, samples=70, samples y=70] function { (1 - 0.4 * 25 * exp(1-(x*x+y*y)) / (8 * 1.4 * pi * pi) )^(2.5) };
    \end{axis}
\end{tikzpicture}

%% file: chapters/Results.tex
\section{ADER-DG in lower precision}\label{chapter:results}

We now investigate the impact of the floating-point precision on the solution computed by ExaHyPE's ADER-DG solver for each of the scenarios presented in chapter~\ref{chapter:verification}. 
In section~\ref{subseq:init_error}, we determine a reference level for each scenario by measuring for each precision how the round-off error already affects the initial conditions 
(this means we disregard errors occurring during the computation of the time-dependent ADER-DG solution). 


In section~\ref{subseq:conv}, we investigate the influence of floating-point precision on the $p$-convergence of the method, i.e., with respect to the polynomial order of the method. 
We expect that for lower precision, convergence will trail off (at the determined reference level, at the latest) as floating-point errors will eventually prevent the algorithm from achieving higher accuracy. We also evaluate whether floating-point errors already affect the asymptotic improvement of accuracy.  


In sections~\ref{subsubseq:isentropic_vortex} and \ref{subsubseq:resting_lake}, we analyse the isentropic density vortex and the lake at rest scenario -- two static scenarios that are (in different aspects) sensitive to correct numerical treatment. 
We thus examine whether using low precision affects the ability of ADER-DG to correctly capture such critical scenarios. 


All of these scenarios were compiled using version $13.3.0$ of the \texttt{g++}-compiler provided by the GNU Compiler Collection and tested on a machine with an Intel i7-10700 CPU.

\subsection{Impact of precision on polynomial representations}\label{subseq:init_error}


The initial errors for the elastic planar wave (section~\ref{subsubseq:planar_waves}) and Euler Gaussian bell (section~\ref{subsubseq:gaussian_bell}) scenarios are represented in figures \ref{fig:initial_error_elastic} and \ref{fig:initial_error_euler}, respectively. 
Each of these were computed for polynomial orders 1 to 9, and either 9 or 27 cells per dimension, corresponding to cell sizes of 0.22 and 0.074, respectively. 
We used the \texttt{fp64} solution as reference, and show the L2 error for the \texttt{fp32}, \texttt{fp16} and \texttt{bf16} solution.

For all precisions, the initial errors do not vary substantially for varying cell sizes and polynomial orders: 
increasing the resolution or polynomial order leads to using more integration points and respective point values, the integrated roundoff error  therefore remains nearly the same. 
We can relate the respective error levels of the initial solution to the errors obtained later for the time-dependent simulations, giving us a more accurate separation between roundoff errors caused by casting the analytical solution to a certain precision, and the numerical and roundoff errors caused when advancing the initial conditions in time with the ADER-DG method.
In particular, this means that any simulations performed entirely in lower precision will result at best in an error determined by this reference level.

Between the different precisions, the errors vary strongly, with \texttt{fp32} being several orders of magnitude closer to the reference solution than \texttt{fp16}, which is again about a factor 10 closer to the solution than \texttt{bf16}. 
The error levels closely match the size of the mantissa of each precision, as depicted in table \ref{table:precision_bits}:
\texttt{fp32} has 23 significant bits ($2^{-23} \approx 10^{-7}$), 
\texttt{fp16} has only 10 significant bits ($2^{-10} \approx 10^{-3}$)
and \texttt{bf16} even less, only 7 bits ($2^{-7} \approx 10^{-2}$). 

The results also confirm that the number of exponent bits of \texttt{fp16} is sufficient to represent the analytic solutions.
\texttt{bf16} can represent a wider range of numbers and also values close to 0 with higher accuracy than \texttt{fp16}, but does not profit from these properties in our scenarios, as the planar wave scenario as well as the Gaussian density bell take values well within \texttt{fp16}'s range of representation.

\input{figures/initial_errors/elastic/elastic_initial}
\input{figures/initial_errors/euler/euler_initial}

\subsection{Impact of precision on convergence}\label{subseq:conv}

We now evaluate the impact of precision on the convergence of the solution for the acoustic and elastic planar wave scenarios (section~\ref{subsubseq:planar_waves}), as well as for the Euler Gaussian bell scenario (section~\ref{subsubseq:gaussian_bell}). 
We computed each solution for each precision, for each polynomial order from 1 to 9 ($p$-convergence), and for either 9 or 27 cells per spatial dimension. The results are shown in figures \ref{fig:polynomial_error_acoustic}, \ref{fig:polynomial_error_elastic} and \ref{fig:polynomial_error_euler}.

\input{figures/convergence/acoustic/acoustic}
\input{figures/convergence/elastic/elastic}
\input{figures/convergence/euler/euler}

We observe the ``classical'' convergence behavior only for \texttt{fp64} precision:
the L2 error decreases with the appropriate convergence order, and the finer-mesh discretisation leads to a substantial improvement of the error. The error eventually plateaus: at $\approx10^{-11}$ for the planar wave scenarios, for polynomial order 6 and higher, and already at $\approx10^{-9}$ for the Euler Gaussian bell, starting from polynomial order 5.


For \texttt{fp32}, we observe larger overall errors than for \texttt{fp64}, and convergence does plateau earlier -- for the planar wave scenarios at \textbf{$\approx10^{-5}$} and already from polynomial order 4.
At this point, the numerical precision apparently becomes the leading source of error -- note that compared to the results in figures~\ref{fig:initial_error_elastic} and \ref{fig:initial_error_euler}, we lose about 1--2 orders of magnitude due to rounding errors throughout the computation of the solution.
For the Euler Gaussian bell, errors also no longer improve for polynomial orders bigger than 4. Here, the \texttt{fp64} solution already no longer improves from polynomial order 5. 
Hence for all three scenarios the largest useful polynomial order is lower by 1 or 2 for \texttt{fp32} than for \texttt{fp64}.

For \texttt{fp16} and \texttt{bf16}, we can no longer recognize any kind of convergence behavior. And neither was capable of producing results for all of the evaluated scenarios. The \texttt{bf16} results for the linear scenarios were still ``correct within an order of magnitude'', that is to say the initial sinusoidal wave was still reflected in the results, but its values were off (about $10^{-1}$ for the $27\times 27$ grid). Specifically, we observed strong dissipation of the wave, with amplitudes decreased.
For \texttt{fp16} and order 2, the L2 error was $\approx 4\cdot 10^{-2}$; same as the maximum error over all quadrature points. This is small enough as to be ``visibly correct''. For polynomial orders from 3 upwards however, \texttt{fp16} is incapable of resolving either of the planar wave scenarios, only computing \texttt{NaN} values (cmp.\ missing results in figure~\ref{fig:polynomial_error_acoustic}). 
We therefore recomputed all simulations for \texttt{fp16} and \texttt{bf16} using \texttt{fp32} precision for the predictor step. Mixed-precision results will be discussed in more detail in chapter \ref{chapter:mixed_precision} -- for the two linear scenarios this enabled the algorithm to compute decently accurate results without failure, though it did not restore $p$-convergence. 
As with the solutions for order 2, we observed maximum errors of around $10^{-2}$, which are not distinctly visible. 
Hence, while these errors are orders of magnitude larger than those produced by \texttt{fp32} and \texttt{fp64}, at least a roughly correct solution could be achieved through the use of mixed-precision.

For the Euler scenario, \texttt{bf16} produced \texttt{NaN}-results beyond a polynomial order of 3. As we will see in section \ref{subseq:static}, here rounding errors cause oscillations in the results. These then diverge until the floating-point formats cannot represent their values anymore, leading first to positive and negative infinite values and eventually to \texttt{NaN} values.
\texttt{fp16} produced valid solutions, but does not show any convergence with polynomial order -- the lowest error is achieved for order 2. 
Where both 16-bit formats worked, \texttt{fp16} was more accurate by roughly one digit due to its larger mantissa causing smaller rounding errors.
This is different to the linear wave scenarios, where \texttt{bf16} was more robust than \texttt{fp16}, indicating that the wider range of representable values, particularly values near 0, was required there. 
For the nonlinear Euler Gaussian-bell scenario, however, \texttt{fp16} could compute a solution but \texttt{bf16} failed, which indicates that the values remained within the representable range of \texttt{fp16}, where its longer mantissa allowed it to be more robust.

Increasing the polynomial order or the spatial discretization did not improve the solution for the 16-bit computations. 
Errors rather deteriorated or even diverged; in the best cases the error remained constant. 
This is most likely caused by the higher number of operations required to simulate scenarios with higher polynomial orders (with higher computational effor per element) or smaller cell sizes (and thus more required timesteps due to the CFL condition). 
Roundoff errors accumulate and increase with increasing number of operations in the ADER-DG kernels, eventually exceeding the respective discretization error. 


\subsection{Static scenarios}\label{subseq:static}

For the Euler isentropic vortex (section~\ref{subsubseq:isentropic_vortex}) and the shallow water lake at rest scenario (section~\ref{subsubseq:resting_lake}), the analytical solution remains constant, as at any point in the domain fluxes cancel out exactly. 
Numerical or floating-point errors, however, might introduce slight distortions in the solution, which may increase over time.
We simulated both scenarios on a very coarse grid of $9\times 9$ cells and with polynomial order~5. 

\begin{figure}
    \begin{tabular}{cc}
        \includegraphics[width=0.5\columnwidth]{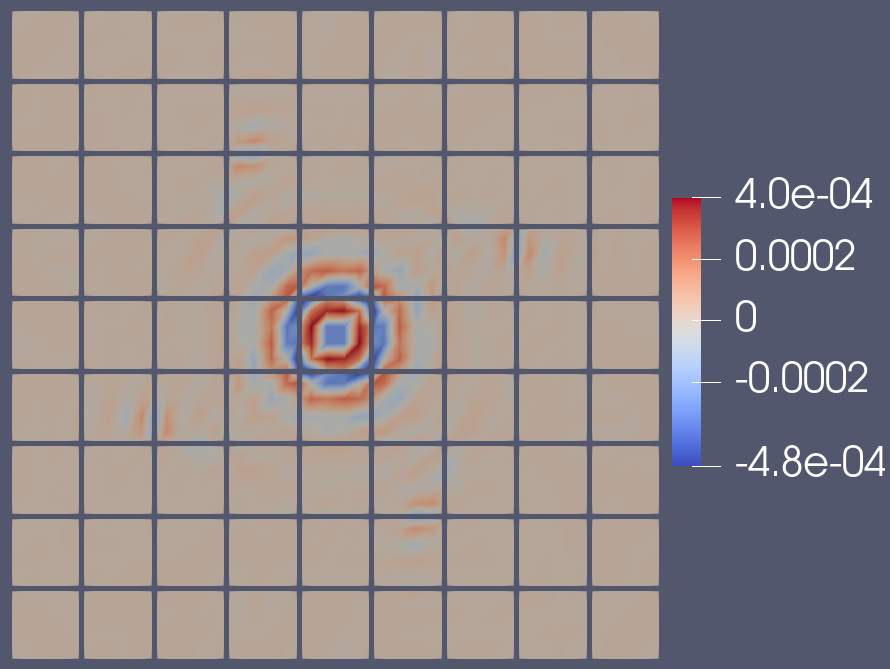} &
        \includegraphics[width=0.5\columnwidth]{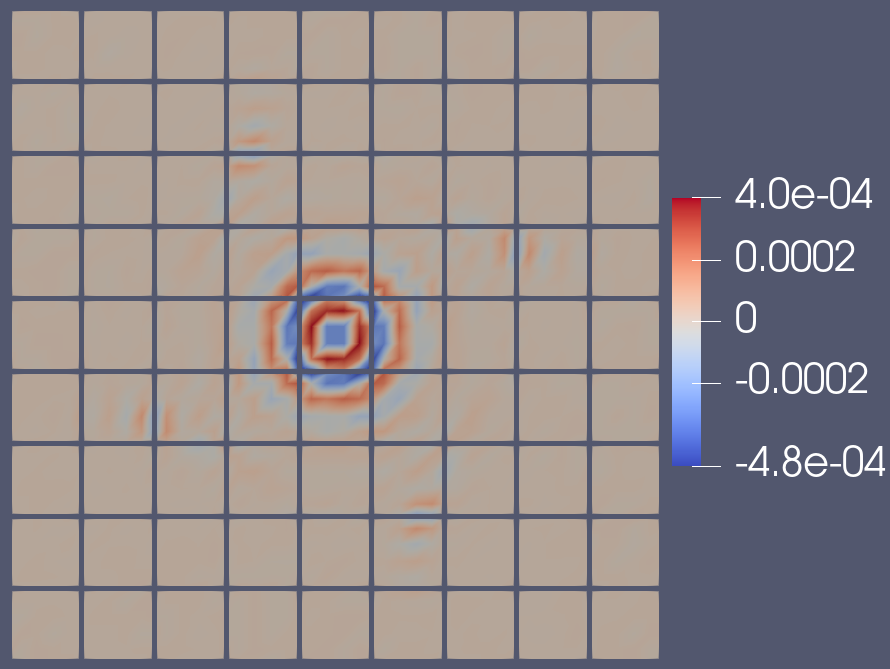} \\
        \textbf{fp64} & \textbf{fp32} \\[1em]
        \includegraphics[width=0.5\columnwidth]{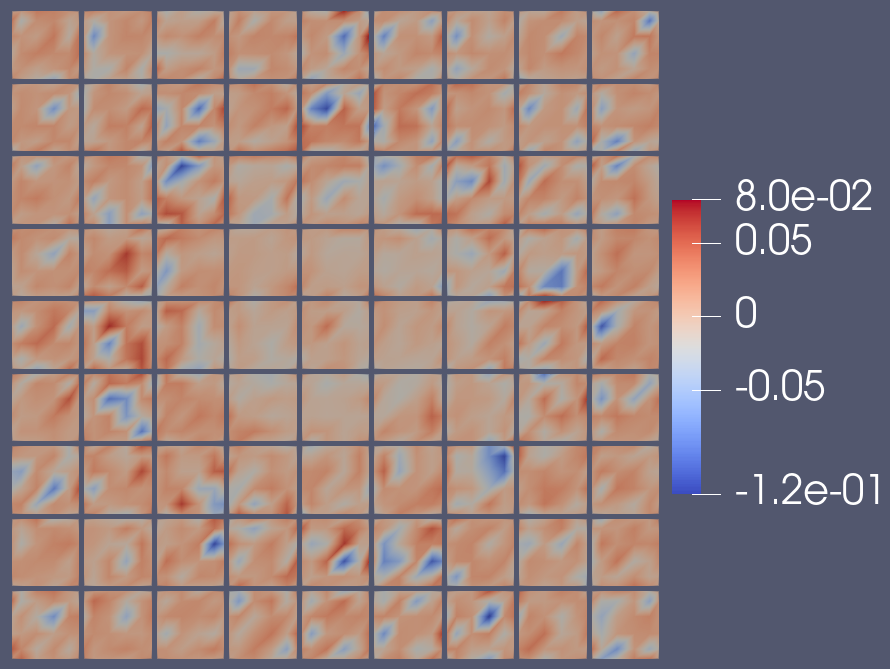} &
        \includegraphics[width=0.5\columnwidth]{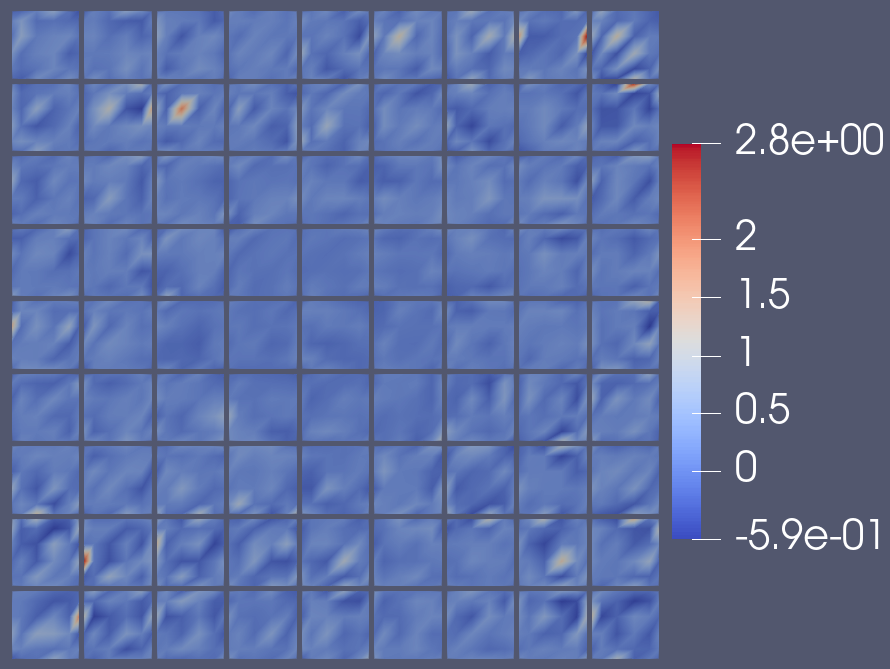} \\
        \textbf{fp16} & \textbf{bf16} \\
    \end{tabular}

    \caption{\label{fig:static_euler_l2}Euler isentropic vortex: 
    Difference of the entropy to that of the analytical solution (at $t=10$; $9\times 9$ grid cells, cf.~the black lines; polynomial order 5). 
    In \texttt{fp64}- and \texttt{fp32}-precision, almost identical errors appear, especially in the vortex centre, indicating that discretization errors dominate. 
    In \texttt{fp16} and \texttt{bf16},  the center of the vortex is still recognisable, but large errors appear over the entire domain -- indicating that floating-point errors are the main source of error.}
\end{figure}

\paragraph{Euler isentropic vortex}
For the isentropic vortex scenario (figure \ref{fig:static_euler_l2}), the errors in \texttt{fp64} and \texttt{fp32} are nearly identical, and focus around the vortex itself. These are apparently caused by discretisation errors as the projected numerical fluxes from the DG discretisation don't cancel out perfectly. 
In contrast, \texttt{fp16} and \texttt{bf16} both produced large errors that show no visible correlation with the initial vortex. 
The errors appear over the entire domain and increase over time, causing the edges of the vortex to become unsteady and rotate with the circular flow. 
Here, the ADER-DG algorithm is apparently incapable of solving the underlying Euler equations due to limited numerical precision. This is also consistent with results for the Euler Gaussian bell scenario (figure~\ref{fig:polynomial_error_euler}), which \texttt{bf16} could not compute and for which \texttt{fp16} exhibited no convergence.

\begin{figure}
    \begin{tabular}{cc}
        \includegraphics[width=0.5\columnwidth]{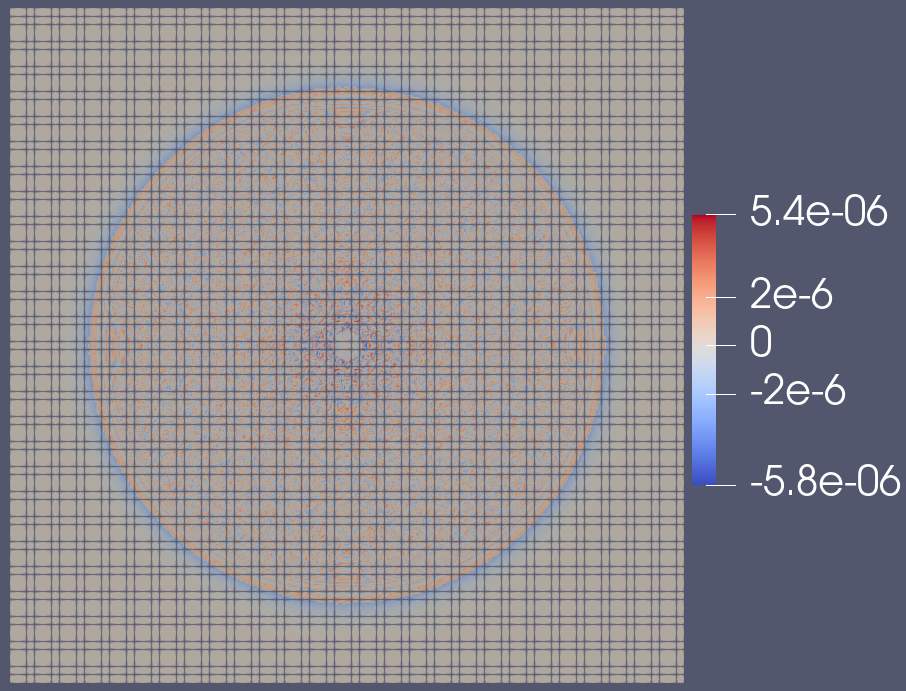} &
        \includegraphics[width=0.5\columnwidth]{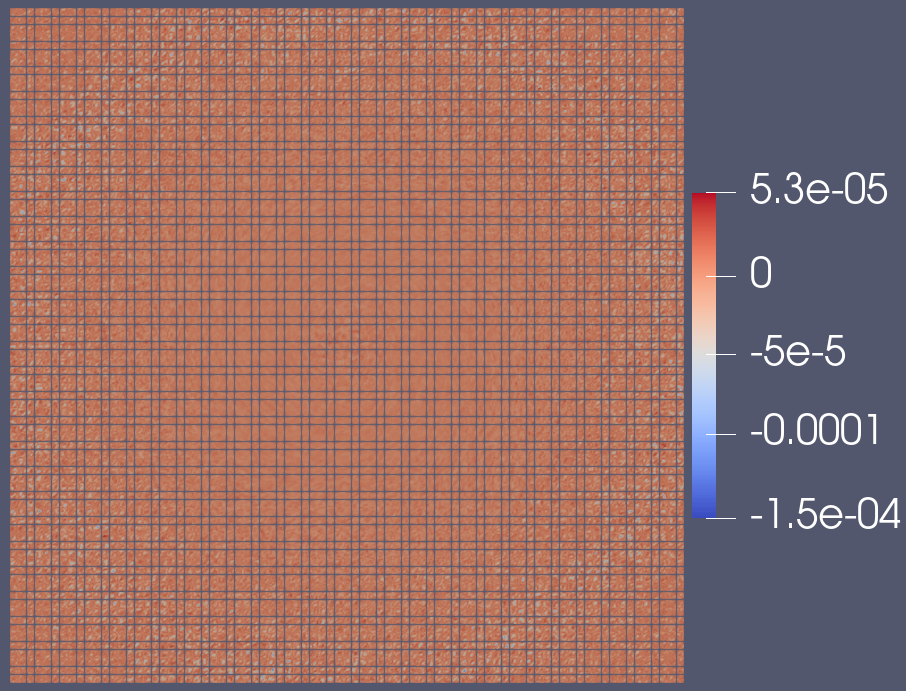} \\
        \textbf{fp64} & \textbf{fp32} \\
    \end{tabular}

    \caption{\label{fig:static_euler_l4}Euler isentropic vortex: 
    Difference of the entropy to that of the analytical solution (at $t=10$; $81\times 81$ grid cells, cf.~the black lines; polynomial order 5). 
    Compared to the results in figure \ref{fig:static_euler_l2}, the now finer resolution leads to much smaller discretization errors. 
    Floating-points errors now become visible for \texttt{fp64}, and dominate the errors for \texttt{fp32}-precision. }
\end{figure}

To further examine the impact of discretization vs.\ floating point errors (for \texttt{fp64} and \texttt{fp32}), we re-computed the isentropic vortex on a finer mesh ($81\times 81$) -- see figure \ref{fig:static_euler_l4}. 
We now observe lower overall errors for both precisions, but also different behavior for \texttt{fp64} and \texttt{fp32}. 
In \texttt{fp64}, the vortex structure is still visible, but random patterns also occur, indicating floating point errors. 
For \texttt{fp32}, the floating point errors seem to be dominant. 
The results confirm that the rounding errors become dominant once the discretisation errors have decreased sufficiently.

The behaviour of the error over time (not representable in the plots) is similar for each precision, and shows a steady increase of the error over the entire duration of the simulation. The increase rate differs for each precision but remains roughly constant over the entire simulation, also for cases where it eventually causes failure. 
This again indicates that errors in the results are caused by accumulation of rounding errors, which -- while different for each precision -- remain within the same order of magnitude within a given precision.

\paragraph{Shallow water lake at rest}

\begin{figure}
    \begin{tabular}{cc}
        \includegraphics[width=0.5\columnwidth]{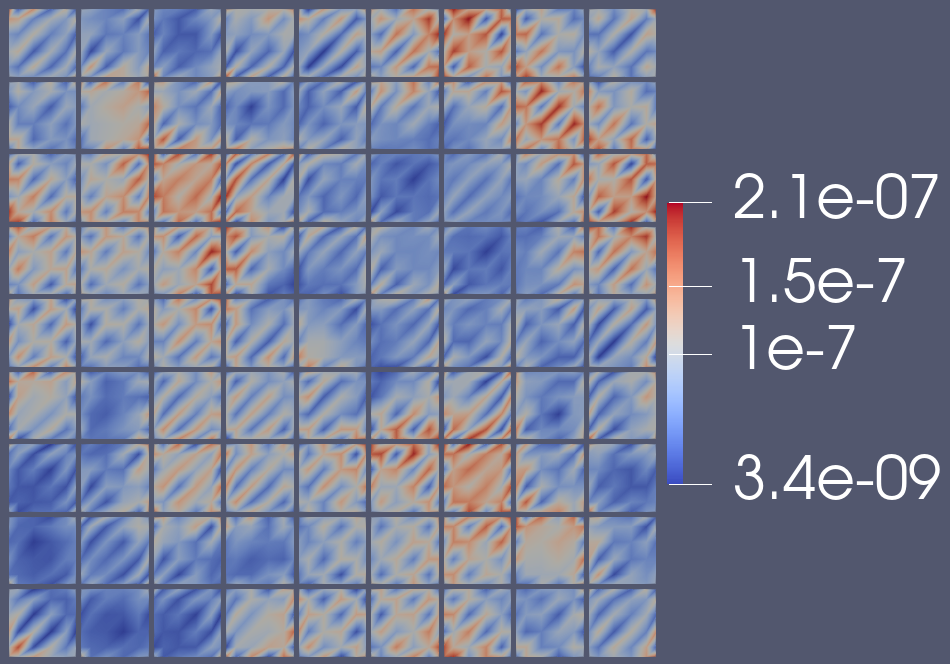} &
        \includegraphics[width=0.5\columnwidth]{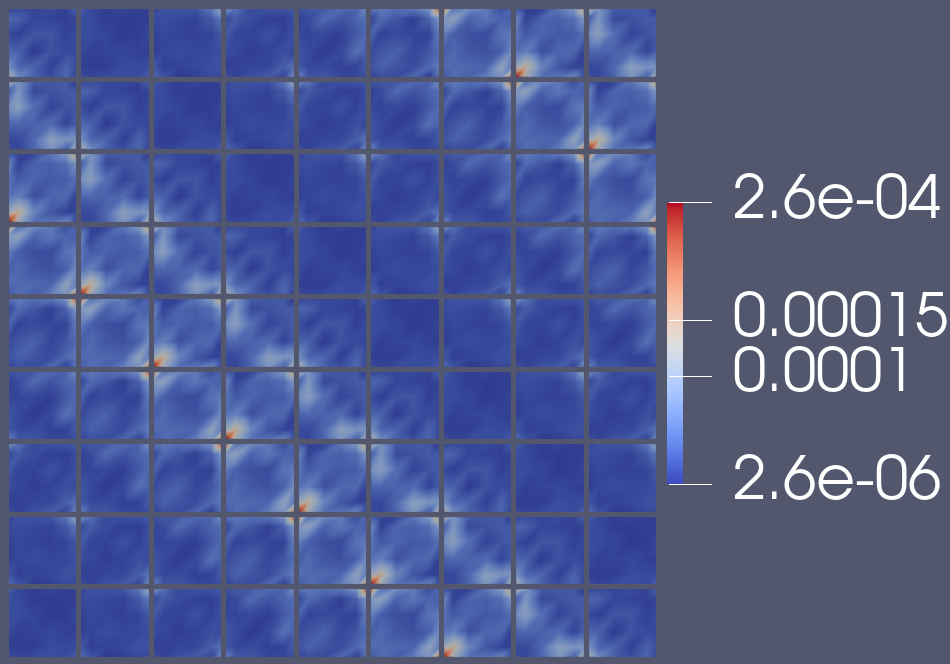} \\
        \textbf{fp64} & \textbf{fp32} \\

        \includegraphics[width=0.5\columnwidth]{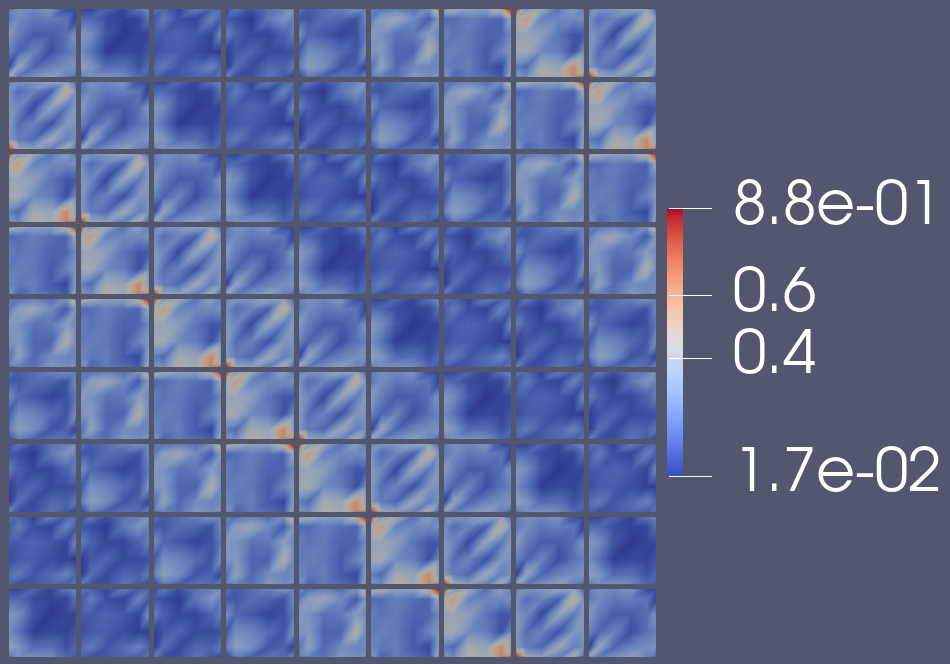} &
        \includegraphics[width=0.5\columnwidth]{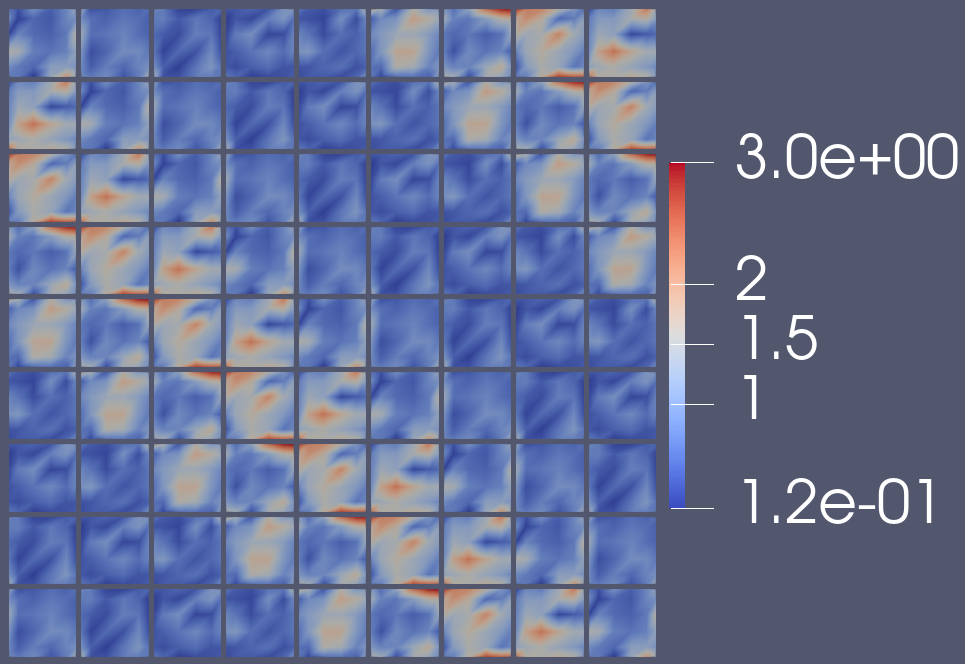} \\
        \textbf{fp16} & \textbf{bf16} \\
    \end{tabular}

    \caption{\label{fig:static_SWE}
    Absolute velocities, at $t=1$, of the lake at rest scenario with sinusoidal bathymetry 
    ($9\times 9$ grid cells, cf.~the black lines; polynomial order 5). 
    In all precisions spurious velocities appear, especially around the crest of the sine function. 
    Within cells, errors are highest when the crest of the sine function passes near the corner of the cell 
    (Runge oscillations, consisting of high-order polynomials with small but near-random coefficients caused by rounding errors). 
    Spurious velocities are small for \texttt{fp64} ($\sim 10^{-7}$) and tolerable for \texttt{fp32} ($\sim 10^{-4}$), but reach unacceptable levels for \texttt{fp16} ($\sim 0.9$) and \texttt{bf16} ($\sim 3.0$). }
\end{figure}

For the lake at rest scenario, we observe (in figure \ref{fig:static_SWE}) for all precisions that spurious velocities and elevations of the sea surface appear, following the profile of the sinusoidal bathymetry $b$.
The errors are largest in the areas with largest gradient, which suggests that rounding errors prevent the well-balanced property (cmp.\ discussion in section~\ref{sec:swe}) from being satisfied exactly, e.g., the terms $\frac{\partial (b + h)}{\partial x}$ and $\frac{\partial (b + h)}{\partial y}$ in equation \eqref{eq:swe} do not sum up exactly to $0$.
Within individual cells, the errors exhibit a pattern typical for the well-known Runge-phenomenon: small deviations in the Lagrange-polynomial oefficients from $0$ lead to oscillations, most strongly near the corners of the cells. 
While the spurious waves exhibit the same pattern, their amplitudes heavily depend on numerical precision, being of the order of $10^{-7}$ in \texttt{fp64}, $10^{-4}$ in \texttt{fp32}, $10^{-1}$ in \texttt{fp16}, and $10^0$ in \texttt{bf16}. 
Hence, induced oscillations are proportional to the chosen precision.

The development of the error over time is again interesting:
For each precision, the error increases initially, until reaching some threshold, where the nearly constant error level depends on the precision. 
Apparently, while spurious wave are created from the initial condition, the simulation approaches a different steady state, which depends on the rounding errors. 

\section{ADER-DG in mixed precision}\label{chapter:mixed_precision}

We now examine where the ADER-DG algorithm can benefit from performing selected kernels in lower or higher precision. 
For this, we simulated all previously presented scenarios (using a 27$\times$27 grid and polynomial order 5) in a prescribed \textbf{storage} precision (\texttt{fp64}, \texttt{fp32}, \texttt{fp16} or \texttt{bf16}; cmp.\ chapter~\ref{chapter:implementation}) but use a higher or lower precision for exactly \emph{one} of the following kernels:
\begin{description}
    \item[predictor] corresponds to the \textbf{fusedSpaceTimePredictorVolumeIntegral()} kernel (cmp.~algorithm~\ref{alg:ader} and chapter~\ref{chapter:implementation}), which combines the element-local space-time predictor and the projection of the predictor values to the cell faces (cf.\ section~\ref{subseq:predictor}). 
    The predictor is a compute-bound kernel and dominates the computing time especially for non-linear PDE systems; computing the predictor in low precision promises speedup due to longer effective SIMD lengths, e.g.
    \item[Picard] iterations are used for nonlinear PDE systems, as a costly substep in the predictor. Same as for the predictor, it is a compute-bound step that can be accelerated through exploiting longer effective SIMD lengths. 
    \item[corrector] updates the solution in each cell, but also computes the Riemann flux on each face and the surface integral over the cell's boundary (cmp.~section~\ref{subseq:corrector} and algorithm~\ref{alg:ader}). 
        The corrector typically does not dominate the runtime of the algorithm, but reducing the precision of face data can lead to performance improvement if the data exchange between neighboring DG elements is a performance-limiting factor. Alternatively, certain equations or scenarios, such as the shallow water equations, require accurate Riemann solvers. Here we might hope to achieve more accurate results at low extra cost, by only increasing the precision of the corrector.
\end{description}
In addition, we examine using higher or lower precision for 
\begin{description}
    \item[all] kernels, i.e., using uniform precision for all kernels, but specifying a separate \textbf{storage} precision --
        here our rationale is that high precision is only necessary to avoid floating point errors during the computation of element updates, but quantities can be stored persistently in lower precision;
        we can hope for speedup if the access to main memory (streaming degrees of freedom into cache and registers) is the performance-limiting factor. This is particularly the case for lower-order scenarios. 
\end{description}

%

\paragraph{Impact of computing kernels in lower precision}
Tables \ref{tab:mixed_prec_storage_fp32} and \ref{tab:mixed_prec_storage_fp64} present for all five scenarios the L2 errors obtained for \texttt{fp32} or \texttt{fp64} as storage precision (and as default precision for all kernels), but using lower precision (\texttt{bf16}, \texttt{fp16}, \texttt{fp32}) for selected kernels.
In line with the results from chapter~\ref{chapter:results}, reducing the precision of the predictor, Picard or corrector kernel always increases the resulting errors. 

Both \texttt{bf16} and \texttt{fp16} suffered from catastrophic failures leading to \texttt{NaN} values, though in different instances. For the linear scenarios, \texttt{bf16} is more robust than \texttt{fp16} and capable of producing results where \texttt{fp16} crashes. 
Recalling the discussion in \ref{subseq:init_error} and table \ref{table:precision_bits}, this is likely due to the lower number of exponent bits of \texttt{fp16}, and rounding of small numbers to zero. Inverting these zero values leads to \texttt{NAN} values.
For the Euler equations, with reduced-precision corrector in particular, the opposite effect appears, with \texttt{fp16} producing valid results where \texttt{bf16} could not, likely because the additional precision of \texttt{fp16} prevented the solution from oscillating out of control.


Regarding the accuracy of the results, \texttt{bf16} is the least accurate, followed by \texttt{fp16}, then \texttt{fp32} -- reflecting the number of mantissa bits. 
Interestingly, reducing the corrector precision consistently had the smallest impact, less severe than predictor and Picard iterations. 
Computing all kernels in reduced precision also consistently increased the error, though it remained within the magnitude of the largest error among individual kernels. 
Occasionally, mixed-precision results using \texttt{fp64} for storage produced worse results than using \texttt{fp32} for storage: 
for example in the acoustic scenario computing the corrector in \texttt{bf16} affected the results in \texttt{fp32} storage precision (L2 error 6.60e{-3}) less than those for \texttt{fp64} (7.02e{-3}). 
However, as with every other case where \texttt{fp32} as storage precision outperforms \texttt{fp64}, errors are very similar, basically indistinguishable when visualised, and multiple orders of magnitude worse than the reference solutions in \texttt{fp32} and \texttt{fp64}. 
Obviously, the reduced-precision kernel dominates the error. 
As a similar example, the shallow water scenario with the corrector computed in \texttt{bf16} produced seemingly valid results for \texttt{fp32}, but not for \texttt{fp64} as storage precision. Here, the mixed-precision errors are again orders of magnitude worse than the reference solution (in \texttt{fp32}, 1.19e{0} vs.\ 8.12e{-5}). 
In fact, the error is so large that the scenario cannot reasonably be considered to have been successfully simulated. 
Visualizing the solution confirms this, as the solution shows large spurious waves. 
The well-balanced property clearly could not be maintained in either case. 

In cases where we do obtain qualitatively correct solutions, computing even the least impactful kernel in a lower precision causes the overall error to increase by several orders of magnitude. 
For example, in the acoustic scenario computing the predictor in \texttt{bf16} led to an increase of the error from 1.07e{-5} to 3.08e{-1} (with \texttt{fp32} as base precision) and from 7.58e{-10} to 3.06e{-1} (with \texttt{fp64}). 
The mixed-precision solutions are qualitatively similar, with the sinusoidal wave being reflected in the results, but noticeably dissipated when compared to the initial conditions. 
The low precision for the predictor obviously dominates the overall error. 
We do not observe any situation where reduced precision for a kernel comes at no cost to the solution.
In some cases, computing one kernel in low precision even results in a higher error than computing the entire algorithm in the low precision. 
Such was the case for the Euler Gaussian bell scenario with Picard iterations in \texttt{fp32} and everything else in \texttt{fp64} (error of 2.50e{-5} vs.\ 2.05e{-5} for everything in \texttt{fp32}). 
The corresponding solutions are hardly distinguishable from another nor from the initial conditions. 
We assume that conversions between precisions cause aspects of the numerics to not quite match. In the static scenarios for example, updates resulting from the predictor and from the corrector must exactly cancel, but might no longer match if different rounding is applied.

\input{figures/mixed_precision/table_mixed_fp32storage}

\input{figures/mixed_precision/table_mixed_fp64storage}

\paragraph{Impact of computing kernels in higher precision}
Next, we compute each scenario with a low default precision, either \texttt{bf16}, \texttt{fp16} or \texttt{fp32}, but switch one of the kernels to \texttt{fp64} precision. The results are presented in table \ref{tab:mixed_prec_higher_fp64}.

We first examine where using high precision in selected kernels can ``heal'' a scenario for which low precision fails. 
There are five instances, in which a scenario cannot be computed with half precision: 
the acoustic, elastic and shallow water scenarios lead to \texttt{NaN} in \texttt{fp16}, and both Euler scenarios fail in \texttt{bf16}. 
For the Euler scenarios and \texttt{bf16}, we observe that switching to double precision for any or even all of the kernels does not help. Apparently, \texttt{bf16} simply does not have a sufficient precision to resolve these scenarios. As the \texttt{fp16} computations succeed, we postulate the lack of mantissa bits in \texttt{bf16} as the main reason. 
The acoustic and elastic planar wave scenarios, with \texttt{fp16} as base precision, can be cured by using higher precision in the predictor. Using higher precision only in the corrector is not sufficient. 
As already indicated in section \ref{subseq:conv}, the mixed-precision \texttt{fp16} solutions for the acoustic and elastic equations are ``correct within an order of magnitude'' (at 9.01e-2 and 2.21e-1, respectively), meaning that the sinusoidal wave is still reflected in the results but noticeably dissipated. 
Using high precision in predictor and corrector leads only to small further improvements. 
We infer that, for the wave equations, a larger range of values is necessary for computing the predictor, which explains why \texttt{bf16} is able to produce solutions (though inaccurate) for acoustic and elastic.  

For the shallow water equations, it is necessary to compute all kernels in increased precision to ``heal'' the \texttt{fp16} computation. Neither increasing the precision for the predictor only, nor for the corrector only, is sufficient. Again the mixed-precision solution can be considered as ``correct within an order of magnitude'', with an L2 error of 4.70e-2. 
For the shallow water equation in \texttt{bf16}, we observed the counter-intuitive situation that increasing the precision for the kernels \texttt{Picard} and \texttt{corrector} causes failure. 
However, as shown in figure \ref{fig:static_SWE}, the errors in \texttt{bf16} were so large (L2 error at 1.07e0) that the scenario cannot be considered to be successfully computable. 
We therefore consider these \texttt{NaN} failures rather as a failure to ``heal'' the scenario with mixed precision than as causing a successful scenario to fail. 
Recalling again previous results from table \ref{tab:mixed_prec_storage_fp64}, that computing any kernel in \texttt{bf16} causes a computation with \texttt{fp64} as storage precision to fail, we see this as confirmation that \texttt{bf16} is unsuitable to simulate the shallow water lake at rest scenario, due to its delicate numerics. 
Even computing all kernels in \texttt{fp64} resulted in an L2 error of 4.05e{-1} -- with spurious waves larger than $10\%$ of the water height, this is far away from the requirements typically posed by scenarios that rely on a well-balanced scheme. By contrast, using \texttt{fp16} as storage precision while computing all kernels in \texttt{fp64} resulted in an error of 4.70e{-2}. Here the waves are within about $1\%$ of the water height, and the solution closely resembles the initial conditions, up to small spurious oscillations around the crest of the sine.

%
\input{figures/mixed_precision/table_mod_fp64_all}

Next, we examine where using higher precision for selected kernels can lead to improved accuracy of a (successful) low-precision computation. 
Here, we report a clearly negative result for the acoustic and elastic wave equation: increasing the precision of any given kernel does not improve the resulting solution. Typically the L2 error remained within the same magnitude of the base error, regardless of the base precision (\texttt{bf16}, \texttt{fp16} or \texttt{fp32}).
This is consistent with the fact that all kernels of the linear PDE solvers are essentially linear operators, and also do not deviate too much in terms of computational effort. Hence, there is no kernel that strongly dominates the accuracy of the final result.
The strongest impact comes from storage precision. 

For the non-linear PDEs, the predictor with its costly Picard loop strongly dominates the computational effort, and we have scenarios where modified Riemann solvers are required -- giving us clear ``suspects'' for kernels that might have more impact on the error than others. 
However, also for the non-linear scenarios, improvements due to mixed-precision are limited. 
In the shallow water equations, we apply a custom Riemann solver \eqref{eq:swe_custom_riemann} to ensure well-balancedness. 
Our results show that we get a substantial improvement in the L2 error (approx.~one order of magnitude) only when computing all kernels in high precision. The resulting error level also stays one order of magnitude below the \texttt{bf64} error, though. 
Here it could be interesting to see the respective impact when using more complicated Riemann solvers, such as those by George et al.~\cite{GeorgeRiemann}. 
Finally, for the two Euler scenarios, we also get mixed impressions: 
for the Gaussian bell scenario, we gain roughly one order of magnitude for \texttt{fp32}, when computing the Picard loop or all the kernels in high precision. We do not see such an improvement for the isentropic vortex. 
For table~\ref{tab:mixed_prec_higher_fp64_depth_4}, we recomputed the Euler scenarios on a finer mesh, to rule out that errors are limited by the discretization for error. 
We again see the mixed-precision gain only for the Gaussian bell scenario. 

\input{figures/mixed_precision/table_mod_fp64_depth_4}

%% file: figures/initial_errors/elastic/elastic_initial.tex
\definecolor{lightGreen}{rgb}{.2,.8,.6}

\begin{figure}[ht]
    \centering
    \scalebox{0.50}{\begin{tikzpicture}[font=\LARGE]
        \begin{axis}[
            xlabel=polynomial order,
            ylabel=l2 error,
            ymin=0.000000335, ymax=0.000000365,
            enlargelimits = true,
            xtick=data
           ]
            \addplot[purple,thick,mark=otimes*] table [y=fp32,x=Order]{figures/initial_errors/elastic/elastic_depth_2.dat};
            \addlegendentry{fp32, 9x9 cells}
            \addplot[violet,thick,mark=otimes*] table [y=fp32,x=Order]{figures/initial_errors/elastic/elastic_depth_3.dat};
            \addlegendentry{fp32, 27x27 cells}]
        \end{axis}
    \end{tikzpicture}}
    \quad
    \scalebox{0.50}{\begin{tikzpicture}[font=\LARGE]
        \begin{axis}[
            xlabel=polynomial order,
            ylabel=l2 error,
            ymin=0.0027, ymax=0.00295,
            enlargelimits = true,
            xtick=data,
            legend pos= north west
           ]
            \addplot[olive,thick,mark=otimes*] table [y=fp16,x=Order]{figures/initial_errors/elastic/elastic_depth_2.dat};
            \addlegendentry{fp16, 9x9 cells}
            \addplot[brown,thick,mark=otimes*] table [y=fp16,x=Order]{figures/initial_errors/elastic/elastic_depth_3.dat};
            \addlegendentry{fp16, 27x27 cells}]
        \end{axis}
    \end{tikzpicture}}
    \quad
    \scalebox{0.50}{\begin{tikzpicture}[font=\LARGE]
        \begin{axis}[
            xlabel=polynomial order,
            ylabel=l2 error,
            ymin=0.021, ymax=0.026,
            enlargelimits = true,
            xtick=data,
            legend pos= north east
           ]
            \addplot[cyan,thick,mark=otimes*] table [y=bf16,x=Order]{figures/initial_errors/elastic/elastic_depth_2.dat};
            \addlegendentry{bf16, 9x9 cells}
            \addplot[blue,thick,mark=otimes*] table [y=bf16,x=Order]{figures/initial_errors/elastic/elastic_depth_3.dat};
            \addlegendentry{bf16, 27x27 cells}]
        \end{axis}
    \end{tikzpicture}}

    \caption{L2 error for the initial solution of the elastic planar-wave scenario (Sec.~\ref{subsubseq:planar_waves}) with different polynomial orders, cell sizes and precisions.}
    \label{fig:initial_error_elastic}
\end{figure}

%% file: figures/initial_errors/euler/euler_initial.tex
\definecolor{lightGreen}{rgb}{.2,.8,.6}

\begin{figure}[htb]
    \centering
    \scalebox{.50}{\begin{tikzpicture}[font=\LARGE]
        \begin{axis}[
            xlabel=polynomial order,
            ylabel=l2 error,
            ymin=0.0000000625, ymax=0.0000000690,
            enlargelimits = true,
            xtick=data
           ]
            \addplot[purple,thick,mark=otimes*] table [y=fp32,x=Order]{figures/initial_errors/euler/euler_depth_2.dat};\addlegendentry{fp32, 9x9 cells}
            \addplot[violet,thick,mark=otimes*] table [y=fp32,x=Order]{figures/initial_errors/euler/euler_depth_3.dat};
            \addlegendentry{fp32, 27x27 cells}]
        \end{axis}
    \end{tikzpicture}}    
    \quad
    \scalebox{0.50}{\begin{tikzpicture}[font=\LARGE]
        \begin{axis}[
            xlabel=polynomial order,
            ylabel=l2 error,
            ymin=0.000525, ymax=0.00055,
            enlargelimits = true,
            xtick=data,
            legend pos= north west
           ]
            \addplot[olive,thick,mark=otimes*] table [y=fp16,x=Order]{figures/initial_errors/euler/euler_depth_2.dat};
            \addlegendentry{fp16, 9x9 cells}
            \addplot[brown,thick,mark=otimes*] table [y=fp16,x=Order]{figures/initial_errors/euler/euler_depth_3.dat};
            \addlegendentry{fp16, 27x27 cells}]
        \end{axis}
    \end{tikzpicture}}
    \quad
    \scalebox{0.50}{\begin{tikzpicture}[font=\LARGE]
        \begin{axis}[
            xlabel=polynomial order,
            ylabel=l2 error,
            ymin=0.00468, ymax=0.00471,
            enlargelimits = true,
            xtick=data,
            legend pos= north east
           ]
            \addplot[cyan,thick,mark=otimes*] table [y=bf16,x=Order]{figures/initial_errors/euler/euler_depth_2.dat};
            \addlegendentry{bf16, 9x9 cells}
            \addplot[blue,thick,mark=otimes*] table [y=bf16,x=Order]{figures/initial_errors/euler/euler_depth_3.dat};
            \addlegendentry{bf16, 27x27 cells}]
        \end{axis}
    \end{tikzpicture}}

    \caption{L2 error for the initial solution of the Euler Gaussian bell scenario (Sec.~\ref{subsubseq:gaussian_bell}) with different polynomial orders, cell sizes and precisions.}
    \label{fig:initial_error_euler}
\end{figure}

%% file: figures/convergence/acoustic/acoustic.tex
\definecolor{lightGreen}{rgb}{.2,.8,.6}

\begin{figure}
    \centering
    \scalebox{0.50}{\begin{tikzpicture}[font=\LARGE]
        \begin{semilogyaxis}[
            xlabel=polynomial order,
            ylabel=l2 error,
            enlargelimits = true,
            xtick=data
           ]
            \addplot[lightGreen,thick,mark=otimes*] table [y=fp64,x=Order]{figures/convergence/acoustic/acoustic_depth_2.dat};
            \addlegendentry{fp64, 9x9 cells}
            \addplot[teal,thick,mark=otimes*] table [y=fp64,x=Order]{figures/convergence/acoustic/acoustic_depth_3.dat};
            \addlegendentry{fp64, 27x27 cells}
        \end{semilogyaxis}
    \end{tikzpicture}}
    \quad
    \scalebox{0.50}{\begin{tikzpicture}[font=\LARGE]
        \begin{semilogyaxis}[
            xlabel=polynomial order,
            ylabel=l2 error,
            enlargelimits = true,
            xtick=data
           ]
            \addplot[purple,thick,mark=otimes*] table [y=fp32,x=Order]{figures/convergence/acoustic/acoustic_depth_2.dat};
            \addlegendentry{fp32, 9x9 cells}
            \addplot[violet,thick,mark=otimes*] table [y=fp32,x=Order]{figures/convergence/acoustic/acoustic_depth_3.dat};
            \addlegendentry{fp32, 27x27 cells}]
        \end{semilogyaxis}
    \end{tikzpicture}}
    \quad
    \scalebox{0.50}{\begin{tikzpicture}[font=\LARGE]
        \begin{semilogyaxis}[
            xlabel=polynomial order,
            ylabel=l2 error,
            ymin=0.01, ymax=9,
            xtick={1,2,3,4,5,6,7,8,9},
            enlargelimits = true,
            legend pos= north west
           ]
            \addplot[olive,thick,mark=otimes*] table [y=fp16,x=Order]{figures/convergence/acoustic/acoustic_depth_2.dat};
            \addlegendentry{fp16, 9x9 cells}
            \addplot[brown,thick,mark=otimes*] table [y=fp16,x=Order]{figures/convergence/acoustic/acoustic_depth_3.dat};
            \addlegendentry{fp16, 27x27 cells}]
            \addplot[white] coordinates {
                (1,0.0395) (2,0.0395) (3,0.0395) (4,0.0395) (5,0.0395) (6,0.0395) (7,0.0395) (8,0.0395) (9,0.0395)
            };
        \end{semilogyaxis}
    \end{tikzpicture}}
    \quad
    \scalebox{0.50}{\begin{tikzpicture}[font=\LARGE]
        \begin{semilogyaxis}[
            xlabel=polynomial order,
            ylabel=l2 error,
            enlargelimits = true,
            xtick=data,
            legend pos= north west,
            ymin=0.01, ymax=9,
           ]
            \addplot[cyan,thick,mark=otimes*] table [y=bf16,x=Order]{figures/convergence/acoustic/acoustic_depth_2.dat};
            \addlegendentry{bf16, 9x9 cells}
            \addplot[blue,thick,mark=otimes*] table [y=bf16,x=Order]{figures/convergence/acoustic/acoustic_depth_3.dat};
            \addlegendentry{bf16, 27x27 cells}]
        \end{semilogyaxis}
    \end{tikzpicture}}
    \quad
    \scalebox{0.50}{\begin{tikzpicture}[font=\LARGE]
        \begin{semilogyaxis}[
            xlabel=polynomial order,
            ylabel=l2 error,
            enlargelimits = true,
            ymin=0.01, ymax=9,
            xtick=data,
            legend pos= north west
           ]
            \addplot[olive,thick,mark=otimes*] table [y=fp16_mixed,x=Order]{figures/convergence/acoustic/acoustic_depth_2.dat};
            \addlegendentry{\Large{fp16 mixed, 9x9 cells}}
            \addplot[brown,thick,mark=otimes*] table [y=fp16_mixed,x=Order]{figures/convergence/acoustic/acoustic_depth_3.dat};
            \addlegendentry{\Large{fp16 mixed, 27x27 cells}}]
        \end{semilogyaxis}
    \end{tikzpicture}}
    \quad
    \scalebox{0.50}{\begin{tikzpicture}[font=\LARGE]
        \begin{semilogyaxis}[
            xlabel=polynomial order,
            ylabel=l2 error,
            enlargelimits = true,
            ymin=0.01, ymax=9,
            xtick=data,
            legend pos= north west
           ]
            \addplot[cyan,thick,mark=otimes*] table [y=bf16_mixed,x=Order]{figures/convergence/acoustic/acoustic_depth_2.dat};
            \addlegendentry{\Large{bf16 mixed, 9x9 cells}}
            \addplot[blue,thick,mark=otimes*] table [y=bf16_mixed,x=Order]{figures/convergence/acoustic/acoustic_depth_3.dat};
            \addlegendentry{\Large{bf16 mixed, 27x27 cells}}]
        \end{semilogyaxis}
    \end{tikzpicture}}
    \caption{L2 error for the acoustic planar-wave scenario (Sec.~\ref{subsubseq:planar_waves}) after the initial wave has entirely traversed the domain twice -- computed using the ADER-DG method with different polynomial orders, cell sizes and precisions. 
    The \texttt{fp16} computations failed for polynomial orders larger than 2. 
    \texttt{fp16} and \texttt{bf16} results were also computed using mixed precision, with the predictor computed in \texttt{fp32}. Note the different scales of the plots (esp.\ \texttt{fp32} and \texttt{fp64}).}
    \label{fig:polynomial_error_acoustic}
\end{figure}

%% file: figures/convergence/elastic/elastic.tex
\definecolor{lightGreen}{rgb}{.2,.8,.6}

\begin{figure}
    \centering
    \scalebox{0.50}{\begin{tikzpicture}[font=\LARGE]
        \begin{semilogyaxis}[
            xlabel=polynomial order,
            ylabel=l2 error,
            enlargelimits = true,
            xtick=data
           ]
            \addplot[lightGreen,thick,mark=otimes*] table [y=fp64,x=Order]{figures/convergence/elastic/elastic_depth_2.dat};
            \addlegendentry{fp64, 9x9 cells}
            \addplot[teal,thick,mark=otimes*] table [y=fp64,x=Order]{figures/convergence/elastic/elastic_depth_3.dat};
            \addlegendentry{fp64, 27x27 cells}
        \end{semilogyaxis}
    \end{tikzpicture}}
    \quad
    \scalebox{0.50}{\begin{tikzpicture}[font=\LARGE]
        \begin{semilogyaxis}[
            xlabel=polynomial order,
            ylabel=l2 error,
            enlargelimits = true,
            xtick=data
           ]
            \addplot[purple,thick,mark=otimes*] table [y=fp32,x=Order]{figures/convergence/elastic/elastic_depth_2.dat};
            \addlegendentry{fp32, 9x9 cells}
            \addplot[violet,thick,mark=otimes*] table [y=fp32,x=Order]{figures/convergence/elastic/elastic_depth_3.dat};
            \addlegendentry{fp32, 27x27 cells}]
        \end{semilogyaxis}
    \end{tikzpicture}}
    \quad
    \scalebox{0.50}{\begin{tikzpicture}[font=\LARGE]
        \begin{semilogyaxis}[
            xlabel=polynomial order,
            ylabel=l2 error,
            ymin=0.01, ymax=7.5,
            xtick={1,2,3,4,5,6,7,8,9},
            enlargelimits = true,
            legend pos= north west
           ]
            \addplot[olive,thick,mark=otimes*] table [y=fp16,x=Order]{figures/convergence/elastic/elastic_depth_2.dat};
            \addlegendentry{fp16, 9x9 cells}
            \addplot[brown,thick,mark=otimes*] table [y=fp16,x=Order]{figures/convergence/elastic/elastic_depth_3.dat};
            \addlegendentry{fp16, 27x27 cells}]
            \addplot[white] coordinates {
                (1,0.0395) (2,0.0395) (3,0.0395) (4,0.0395) (5,0.0395) (6,0.0395) (7,0.0395) (8,0.0395) (9,0.0395)
            };
        \end{semilogyaxis}
    \end{tikzpicture}}
    \quad
    \scalebox{0.50}{\begin{tikzpicture}[font=\LARGE]
        \begin{semilogyaxis}[
            xlabel=polynomial order,
            ylabel=l2 error,
            enlargelimits = true,
            xtick=data,
            legend pos= north west,
            ymin=0.01, ymax=7.5,
           ]
            \addplot[cyan,thick,mark=otimes*] table [y=bf16,x=Order]{figures/convergence/elastic/elastic_depth_2.dat};
            \addlegendentry{\Large{bf16, 9x9 cells}}
            \addplot[blue,thick,mark=otimes*] table [y=bf16,x=Order]{figures/convergence/elastic/elastic_depth_3.dat};
            \addlegendentry{\Large{bf16, 27x27 cells}}]
        \end{semilogyaxis}
    \end{tikzpicture}}
    \scalebox{0.50}{\begin{tikzpicture}[font=\LARGE]
        \begin{semilogyaxis}[
            xlabel=polynomial order,
            ylabel=l2 error,
            enlargelimits = true,
            ymin=0.01, ymax=7.5,
            xtick=data,
            legend pos= north west
           ]
            \addplot[olive,thick,mark=otimes*] table [y=fp16_mixed,x=Order]{figures/convergence/elastic/elastic_depth_2.dat};
            \addlegendentry{\Large{fp16 mixed, 9x9 cells}}
            \addplot[brown,thick,mark=otimes*] table [y=fp16_mixed,x=Order]{figures/convergence/elastic/elastic_depth_3.dat};
            \addlegendentry{\Large{fp16 mixed, 27x27 cells}}]
        \end{semilogyaxis}
    \end{tikzpicture}}
    \quad
    \scalebox{0.50}{\begin{tikzpicture}[font=\LARGE]
        \begin{semilogyaxis}[
            xlabel=polynomial order,
            ylabel=l2 error,
            enlargelimits = true,
            ymin=0.01, ymax=7.5,
            xtick=data,
            legend pos= north west
           ]
            \addplot[cyan,thick,mark=otimes*] table [y=bf16_mixed,x=Order]{figures/convergence/elastic/elastic_depth_2.dat};
            \addlegendentry{\Large{bf16 mixed, 9x9 cells}}
            \addplot[blue,thick,mark=otimes*] table [y=bf16_mixed,x=Order]{figures/convergence/elastic/elastic_depth_3.dat};
            \addlegendentry{\Large{bf16 mixed, 27x27 cells}}]
        \end{semilogyaxis}
    \end{tikzpicture}}
    \caption{L2 error for the elastic planar-wave scenario (Sec.~\ref{subsubseq:planar_waves}) after the initial wave has entirely traversed the domain twice -- computed using the ADER-DG method with different polynomial orders, cell sizes and precisions. 
    The \texttt{fp16} computations failed for polynomial orders larger than 2. 
    \texttt{fp16} and \texttt{bf16} results were also computed using mixed precision, with the predictor computed in \texttt{fp32}. Note the different scales of the plots (esp.\ \texttt{fp32} and \texttt{fp64}).}
    \label{fig:polynomial_error_elastic}
\end{figure}

%% file: figures/convergence/euler/euler.tex
\definecolor{lightGreen}{rgb}{.2,.8,.6}

\begin{figure}
    \centering
    \scalebox{0.50}{\begin{tikzpicture}[font=\LARGE]
        \begin{semilogyaxis}[
            xlabel=polynomial order,
            ylabel=l2 error,
            enlargelimits = true,
            xtick=data
           ]
            \addplot[lightGreen,thick,mark=otimes*] table [y=fp64,x=Order]{figures/convergence/euler/euler_depth_2.dat};
            \addlegendentry{fp64, 9x9 cells}
            \addplot[teal,thick,mark=otimes*] table [y=fp64,x=Order]{figures/convergence/euler/euler_depth_3.dat};
            \addlegendentry{fp64, 27x27 cells}
        \end{semilogyaxis}
    \end{tikzpicture}}
    \quad
    \scalebox{0.50}{\begin{tikzpicture}[font=\LARGE]
        \begin{semilogyaxis}[
            xlabel=polynomial order,
            ylabel=l2 error,
            enlargelimits = true,
            xtick=data
           ]
            \addplot[purple,thick,mark=otimes*] table [y=fp32,x=Order]{figures/convergence/euler/euler_depth_2.dat};
            \addlegendentry{fp32, 9x9 cells}
            \addplot[violet,thick,mark=otimes*] table [y=fp32,x=Order]{figures/convergence/euler/euler_depth_3.dat};
            \addlegendentry{fp32, 27x27 cells}]
        \end{semilogyaxis}
    \end{tikzpicture}}
    \quad
    \scalebox{0.50}{\begin{tikzpicture}[font=\LARGE]
        \begin{semilogyaxis}[
            xlabel=polynomial order,
            ylabel=l2 error,
            enlargelimits = true,
            ymax=0.18,
            xtick=data,
            legend pos= north west
           ]
            \addplot[olive,thick,mark=otimes*] table [y=fp16,x=Order]{figures/convergence/euler/euler_depth_2.dat};
            \addlegendentry{\Large{fp16, 9x9 cells}}
            \addplot[brown,thick,mark=otimes*] table [y=fp16,x=Order]{figures/convergence/euler/euler_depth_3.dat};
            \addlegendentry{\Large{fp16, 27x27 cells}}]
        \end{semilogyaxis}
    \end{tikzpicture}}
    \caption{L2 error for the Euler Gaussian bell scenario (Sec.~\ref{subsubseq:gaussian_bell}), computed using the ADER-DG method with different polynomial orders, cell sizes and precisions. Note that the scales differ between plots.}
    \label{fig:polynomial_error_euler}
\end{figure}

%% file: figures/mixed_precision/table_mixed_fp32storage.tex
\begin{table}
    \begin{center}
        \begin{tabular}{ |c|| c|c|c| c|c|c| } 
            \hline
                & \multicolumn{3}{|c|}{acoustic planar wave} & \multicolumn{3}{|c|}{elastic planar wave} \\
            \hline
            prec & predictor & corrector & all & predictor & corrector & all \\
            \hline

            bf16 & 3.08e{-1} & 2.02e{-1} & 3.04e{-1} & 4.50e{-1} & 1.85e{-1} & 4.72e{-1} \\

            fp16 & \texttt{NAN} & 6.60e{-3} & \texttt{NAN} & \texttt{NAN} & 1.35e{-2} & \texttt{NAN} \\

            \hline
            fp32 & \multicolumn{3}{|c|}{1.07e{-5}} & \multicolumn{3}{|c|}{1.93e{-5}} \\
            \hline
        \end{tabular}

\medskip
        \begin{tabular}{ |c|| c|c|c|c| c|c|c|c| } 
        \hline
            & \multicolumn{4}{|c|}{Euler Gaussian bell} & \multicolumn{4}{|c|}{Euler isentropic vortex} \\
        \hline
            prec & Picard & predictor & corrector & all & Picard & predictor & corrector & all \\
        \hline

            bf16 & 9.44e{-2} & 1.48e{-1} & \texttt{NAN} & 1.64e-1 & 1.98e0  & 2.50e0 & 2.06e0 & 2.96e0 \\

            fp16 & 2.02e{-2} & 3.22e{-2} & 2.21e{-2} & 3.25e-2 & 1.19e{0} & 2.94e{-1} & 3.35e{-1} & 1.17e0 \\

        \hline
            fp32 & \multicolumn{4}{|c|}{2.05e{-5}} & \multicolumn{4}{|c|}{7.50e{-5}} \\
        \hline
        \end{tabular}
\medskip

        \begin{tabular}{ |c|| c|c|c|c| } 
        \hline
            &
            \multicolumn{4}{|c|}{SWE lake at rest}\\
        \hline
        prec & Picard & predictor & corrector & all \\
        \hline

        bf16 & \texttt{NAN} & \texttt{NAN} & 1.19e{0} & \texttt{NAN} \\

        fp16 & \texttt{NAN} & \texttt{NAN} & \texttt{NAN} & \texttt{NAN} \\

        \hline
        fp32 & \multicolumn{4}{|c|}{8.12e{-5}} \\
        \hline
    \end{tabular}

    \end{center}
    \caption{\label{tab:mixed_prec_storage_fp32}L2 error for all five scenarios, computed with \texttt{fp32} as storage precision, but using lower precision (cf.~column ``prec'') for the Picard kernel, the predictor step, the corrector step or all compute kernels (on a $27\times 27$ grid with polynomial order 5).
    The last row states, as reference, the L2 error for computing with uniform \texttt{fp32} precision.}
\end{table}

%% file: figures/mixed_precision/table_mixed_fp64storage.tex
\begin{table}
    \begin{center}

             \begin{tabular}{ |c|| c|c|c| c|c|c| c|c|c|c| } 
            \hline
                & \multicolumn{3}{|c|}{acoustic} & \multicolumn{3}{|c|}{elastic} \\
            \hline
            prec & predictor & corrector & all & predictor & corrector & all \\
            \hline

            bf16 & 3.06e{-1} & 2.09e{-1} & 2.94e{-1} & 4.50e{-1} & 1.85e{-1} & 4.71e{-1} \\

            fp16 & \texttt{NAN} & 7.02e{-3} & \texttt{NAN} & \texttt{NAN} & 1.35e{-2} & \texttt{NAN} \\
            
            fp32 & 9.83e{-6} & 5.19e{-7} & 9.83e{-6} & 1.58e{-5} & 1.60e{-6} & 1.58e{-5} \\

            \hline
            fp64 & \multicolumn{3}{|c|}{7.58e{-10}} & \multicolumn{3}{|c|}{1.18e{-9}} \\
            \hline
        \end{tabular}

\medskip
        \begin{tabular}{ |c|| c|c|c|c| c|c|c|c| } 
        \hline
            & \multicolumn{4}{|c|}{Euler Gaussian bell} & \multicolumn{4}{|c|}{Euler isentropic vortex} \\
        \hline
            prec & Picard & predictor & corrector & all & Picard & predictor & corrector & all \\
        \hline

            bf16 & 9.34e{-2} & 1.46e{-1} & \texttt{NAN} & 1.63e-1 & 1.99e{0} & 2.50e{0} & 2.08e{0} & 2.94e0 \\

            fp16 & 2.22e{-2} & 3.20e{-2} & 2.28e{-2} & 3.15e-2 & 1.21e{0} & 3.00e{-1} & 3.40e{-1} & 1.16e0 \\

            fp32 & 2.50e{-5} & 2.96e{-6} & 1.68e{-6} & 1.88e{-5} & 3.96e{-5} & 2.91e{-5} & 2.68e{-5} & 4.52e{-5}\\

        \hline
            fp64 & \multicolumn{4}{|c|}{1.03e{-9}} & \multicolumn{4}{|c|}{1.64e{-5}} \\
        \hline
        \end{tabular}
\medskip

    \begin{tabular}{ |c|| c|c|c|c| } 
        \hline
            &
            \multicolumn{4}{|c|}{SWE lake at rest} \\
        \hline
        prec & Picard & predictor & corrector & all \\
        \hline

        bf16 & \texttt{NAN} & \texttt{NAN} & \texttt{NAN} &  1.03e0 \\

        fp16 & \texttt{NAN} & \texttt{NAN} & \texttt{NAN} & \texttt{NAN} \\
        
        fp32 & 3.86e{-5} & 1.36e{-4} & 7.92e{-5} & 8.92e{-5} \\

        \hline
        fp64 & \multicolumn{4}{|c|}{1.63e{-7}} \\
        \hline
    \end{tabular}

    \end{center}
    \caption{\label{tab:mixed_prec_storage_fp64}L2 error for all five scenarios, computed with \texttt{fp64} as storage precision, but using lower precision (cf.~column ``prec'') for the Picard kernel, the predictor step, the corrector step or all compute kernels (on a $27\times 27$ grid with polynomial order 5). 
    The last row states, as reference, the L2 error for computing with uniform \texttt{fp64} precision.}
\end{table}

%% file: figures/mixed_precision/table_mod_fp64_all.tex
\begin{table}
    \begin{center}

        \begin{tabular}{ |c|| c|c|c| c|c|c| } 
            \hline
                & \multicolumn{3}{|c|}{acoustic} & \multicolumn{3}{|c|}{elastic}  \\
            \hline
            base & bf16 & fp16 & fp32 & bf16 & fp16 & fp32 \\
            \hline

            reference & 1.04e{0} & \texttt{NAN} & 1.07e{-5} & 9.41e{-1} & \texttt{NAN} & 1.93e{-5} \\
            \hline


            predictor & 1.20e{0} & 9.01e{-2} & 5.97e{-6} & 7.36e{-1} & 2.21e{-1} & 1.27e{-5} \\

            corrector & 6.56e{-1} & \texttt{NAN} & 1.07e{-5} & 9.59e{-1} & \texttt{NAN} & 1.93e{-5} \\

           all & 7.01e{-1} & 6.70e{-2} & 6.36e{-6} & 7.35e{-1} & 2.12e{-1} & 1.52e{-5} \\

            \hline
        \end{tabular}

\medskip
        \begin{tabular}{ |c|| c|c|c| c|c|c| } 
            \hline
            & \multicolumn{3}{|c|}{Euler Gaussian bell} & \multicolumn{3}{|c|}{Euler isentropic Vortex}\\
            \hline
            base & bf16 & fp16 & fp32 & bf16 & fp16 & fp32 \\
            \hline

            reference & \texttt{NAN} & 2.90e{-2} & 2.05e{-5} & \texttt{NAN} & 1.36e{0} & 7.50e{-5} \\
            \hline

            Picard & \texttt{NAN} &  3.14e{-2} & 3.30e{-6} & \texttt{NAN} & 6.79e{-1} & 6.83e{-5} \\

            predictor & \texttt{NAN} & 5.71e{-2} & 2.19e{-5} & \texttt{NAN} & 1.25e{0} & 7.27e{-5} \\

            corrector & \texttt{NAN} & 2.77e{-2} & 2.05e{-5} & \texttt{NAN} & 1.43e{0} & 7.47e{-5} \\

            all & \texttt{NAN} & 1.98e{-2} & 2.43e{-6} & \texttt{NAN} & 6.83e{-1} & 6.56e{-5} \\

            \hline
        \end{tabular}
\medskip

        \begin{tabular}{ |c|| c|c|c| } 

            \hline
                & \multicolumn{3}{|c|}{SWE lake at rest} \\
            \hline
            base & bf16 & fp16 & fp32 \\
            \hline

            reference & 1.07e{0} & \texttt{NAN} & 8.12e{-5} \\
            \hline

            Picard & \texttt{NAN} & \texttt{NAN}  & 7.99e{-5} \\

            predictor & 9.29e{-1} & \texttt{NAN}  & 6.96e{-5} \\

            corrector & \texttt{NAN} & \texttt{NAN}  & 7.40e{-5} \\

            all & 4.05e{-1} & 4.70e{-2} & 4.12e{-6} \\

            \hline
        \end{tabular}

    \end{center}
    \caption{\label{tab:mixed_prec_higher_fp64}L2 error for all five scenarios, computed using \texttt{fp64} precision for Picard kernel, predictor step, corrector step or persistent storage precision. All other kernels have been performed in the specified ``base'' precision (on a $27\times 27$ grid with polynomial order 5).
    }
\end{table}

%% file: figures/mixed_precision/table_mod_fp64_depth_4.tex
\begin{table}
    \begin{center}

        \begin{tabular}{ |c|| c|c| } 
            \hline
            & Euler Gaussian bell & Euler Isentropic Vortex \\
            \hline
            base & fp32 & fp32 \\
            \hline

            reference &  3.49e-5 & 1.51e-4 \\
            \hline

            Picard    &  3.95-6  & 1.02e-4 \\

            predictor &  3.40e-5 & 1.48e-4 \\

            corrector &  3.41e-5 & 1.50e-4 \\

            all       &  2.78e-6 & 9.82e-5 \\

            \hline
            reference fp64 results & 2.01e-7 & 9.82-6 \\

            \hline
        \end{tabular}
    \end{center}

    \caption{\label{tab:mixed_prec_higher_fp64_depth_4}L2 error for both of the Euler scenarios, computed using \texttt{fp64} precision for Picard kernel, predictor step, corrector step or persistent storage precision. All other kernels have been performed in \texttt{fp32} precision (on a $81\times 81$ grid with polynomial order 5).
    }
\end{table}

%% file: chapters/Conclusion.tex
\section{Conclusion and Findings}\label{chapter:conclusion}

In this work, we evaluated the accuracy achieved by ExaHyPE's ADER-DG implementation, which provides four precisions (\texttt{bf16}, \texttt{fp16}, \texttt{fp32} and \texttt{fp64}) as well as using uniform vs.\ mixed precision. 
We examined four very common hyperbolic PDE systems (acoustic/elastic wave, Euler and shallow water equations) and studied for each the impact of low- and mixed precision by examining the convergence behaviour and the errors obtained for relevant but simple scenarios with known analytical solution. 
From our observations, we summarize the following findings:
\begin{itemize}
    \item To achieve full high-order convergence, \texttt{fp64} precision was required for all kernels. 
    \texttt{fp32} only showed convergence to about order 4 or 5. 
    \texttt{fp16} and \texttt{bf16} usually failed for discretisation order beyond 2. 
    Half precision was not accurate enough for any scenario and failed to achieve proper convergence. 
    \item Precision matters for all kernels. Low precision has a strong impact on the achievable accuracy of the results, regardless of using uniform or mixed precision, and especially at higher orders. 
    \item Catastrophic errors (overflow, \texttt{NaN}) resulting from half precision cannot be isolated to one kernel, but -- depending on the computed scenario -- may occur in any of the kernels. 
    Depending on whether the relative error or range of values is more critical, \texttt{fp16} or \texttt{bf16} were more successful. 
    Still, any use of half precision proved to be risky, even when only used for parts of the computations.
    \item Mixed precision can prevent catastrophic errors by computing individual kernels in higher precision. 
    Where \texttt{fp16} failed, increasing the precision of the predictor step helped most, whereas where \texttt{bf16} failed increasing the storage precision helped most. 
    However, for many half-precision computations, but especially for \texttt{bf16}, results were only qualitatively correct.
    \item Among the ADER-DG kernels, the precision of the corrector generally had the least impact on the resulting solution. However, for the shallow water equations, where the Riemann solver is important for preserving well-balancedness, we observe the opposite effect.
\end{itemize}
We have limited our study to the core ADER-DG algorithm, and have not considered the impact of precision on limiting, or applying lower or higher (mixed) precision adaptively in critical regions. 
We defer these topics to future work, similar as evaluating performance gains due to low and mixed precision, as this highly depends on the hardware utilized. 

To implement mixed-precision ADER-DG in the hyperbolic ExaHyPE PDE engine, we built on its underlying code generation approach for data structures and exploited its template-oriented approach to generate optimized kernels. 
This template-based approach allows the solver code to be generated only in the specified precisions and minimizes the need to maintain multiple variants of the same code.
The implementation is available at \url{https://gitlab.lrz.de/hpcsoftware/Peano/-/tags/2025ExaHyPEMixedPrecisionTOMS}. 
This specific version of the code also contains all details to reproduce the results provided in this work. 
Build specifications are provided in the form of a Docker file to facilitate the installation of the software as well as its required dependencies.